\documentclass[oneside,letterpaper]{amsart}

\usepackage{import}
\usepackage[margin=3cm]{geometry}

\usepackage{graphicx}
\usepackage{tikz,tikz-3dplot}
\usetikzlibrary{calc}

\usepackage{amsthm}

\usepackage[T1]{fontenc}
\usepackage{amsmath}
\usepackage{amssymb}

\usepackage{amsthm}
\usepackage{comment}
\usepackage{enumitem}

\newtheorem{thm}{Theorem}

\newtheorem*{thm*}{Theorem}
\newtheorem{lemma}[thm]{Lemma}
\newtheorem{prop}[thm]{Proposition}
\newtheorem{corr}[thm]{Corollary}
\newtheorem*{corr*}{Corollary}

\theoremstyle{definition}
\newtheorem{defn}{Definition}

\numberwithin{thm}{section}
\numberwithin{defn}{section}

\newcommand{\R}{\mathbb{R}}
\newcommand{\C}{\mathbb{C}}
\newcommand{\of}{\overline{f}}

\newcommand{\E}{\mathcal{E}}

\renewcommand{\Re}{\text{Re}}
\newcommand{\Sch}{Schr\"{o}dinger}
\renewcommand{\Im}{\mathop{\text{Im}}}
\renewcommand{\Re}{\mathop{\text{Re}}}

\newcommand{\gll}{GLL}

\newcommand{\der}[2]{ \frac{ \partial #1 }{ \partial #2 }}
\newcommand{\sder}[2]{ \frac{ \partial^2 #1 }{ \partial #2^2 }}

\newcommand{\Hol}{H\"{o}lder}
\newcommand{\Sp}{\mathbb{S}}
\newcommand{\Cn}[1][n]{\mathbb{C}^{#1}}
\newcommand{\CPn}[1][n]{\mathbb{CP}^{#1}}
\newcommand{\Z}{\mathbb{Z}}

\newcommand{\zkb}{\overline{z_k}}
\newcommand{\At}{\tilde{A}}

\renewcommand{\Re}{\text{Re}\,}
\renewcommand{\Im}{\text{Im}\,}

\newcommand{\Hdot}[1]{\dot{H}^{#1}}

\newcommand{\LebST}[2]{
	L^{#1}_t L^{#2}_x
}

\newcommand{\LebS}[1]{
	{L^{#1}_x}
}

\title[Equivariant Heat and Schr\"{o}dinger flows]{Equivariant Heat and Schr\"{o}dinger flows from Euclidean space to complex projective space}

\author{James Fennell}

\address{Courant Institute of Mathematical Sciences, New York University, 251 Mercer St., New York, NY 11206, United States of America}
\email{jamespfennell@gmail.com}

\begin{document}

\begin{abstract}
We study the equivariant harmonic map heat flow, \Sch{} maps equation, 
and generalized Landau-Lifshitz equation from $\C^n$ to $\CPn$.
By means of a careful geometric analysis, we
	determine a new, highly useful representation of the problem in
        terms of a PDE for radial functions from $\C^n$ to $\Sp^2$.
Using this new representation, we are able to write 
    explicit formulae for the harmonic maps in this context,
    and prove that they all have infinite energy.
We show that the PDEs admit a 
    family of self-similar solutions with smooth profiles; 
    these solutions again have infinite energy, and give an example 
    of regularity breakdown.
Then, using a variant of the Hasimoto transformation applied to our new equation
for the dynamics,
    we prove a small-data global wellposedness result when $n=2$.
This is, to the best of our knowledge, the first global wellposedness
result for \Sch{} maps when the complex dimension of the target 
is greater than one.

In the final section we study a special case of the harmonic map
heat flow corresponding to initial data valued in one great circle.
    We show that the $n=2$ case of this problem is
    a borderline case for the standard classification theory for PDE
    of its type.
\end{abstract}
\maketitle

\tableofcontents

\section{Introduction}


The harmonic map heat flow and the \Sch{} maps equations are natural generalizations of the linear heat and \Sch{} equations where
	the domain and range of the functions considered are manifolds and the Euclidean partial derivatives are replaced by covariant derivatives.
In this article we will be exclusively discussing the setting when the base space is some Euclidean space $\R^d$ and the target is 
	 a K\"{a}hler manifold $N$ with complex structure $J$.
The energy of a map $u:\R^d \rightarrow N$ is defined by the formula,
$
	\E(u) = (1/2) \int_{\R^d} |du|^2 dV.
$
The Euler Lagrange operator $\tau(u)$ corresponding to $\E$ is calculated,
in coordinates, to be $\tau(u) = \sum_{k=1}^d D_k \partial_k u$,
	where the $D_k$ operators are covariant derivatives on $N$.
The \emph{harmonic map heat flow} is then the Cauchy problem given by,
\begin{equation}
	\label{eqn:1:hmhf}
	u_t = \tau(u) =  \sum_k D_k \partial_k u, \;\; u(0) = u_0,
\end{equation}
while the \emph{\Sch{} maps equation} is the Cauchy problem given by,
\begin{equation}
	\label{eqn:1:sf}
	u_t = J \tau(u) =  J \sum_k D_k \partial_k u, \;\; u(0) = u_0.
\end{equation}
One can also consider the \emph{generalized Landau-Lifshitz (GLL) equation},
defined for $\alpha \in [0,\infty)$ and $\beta \in \R$ by,
\begin{equation}
	\label{eqn:1:mf}
	u_t = (\alpha + J) \tau(u) =
        (\alpha + \beta J) \sum_k D_k \partial_k u, \;\; u(0) = u_0;
\end{equation}
this corresponds, when the range is $\C$, to the PDE $u_t = (\alpha+i\beta)\Delta u$.
Let us emphasize that the linearity of the equations in the familiar case when the target is $\C$ is special: in general these problems are nonlinear because of curvature.

The harmonic map heat flow is a well known and extensively studied problem.
It was introduced in \cite{Eells1964} as a tool for studying the existence of harmonic maps.
These are maps which satisfy $\tau(u) = D_k \partial_k u=0$ and correspond to stationary solutions of all of the problems above.
Vast work has been done on the harmonic map heat flow in the subsequent years; see, for example, \cite{Lin2008} for a textbook treatment. 
We mention only that it has been shown that for general $N$ uniqueness of the harmonic map heat flow does not hold,
	and that one way to demonstrate non-uniqueness is through studying self-similar solutions, as is done in \cite{Germain2011,Germain2017}.
This approach is used to prove a non-uniqueness result
for the case of the flow for maps from $\C^2\cong\R^4$ to $\CPn[2]$, 
	 in Section \ref{sec:realheatflow} below.

As opposed to the harmonic map heat flow, the \Sch{} maps equation \eqref{eqn:1:sf} has been much less studied in general.
For the setting we are considering here, that of the flow for maps $u:\R^d \rightarrow N$, local well-posedness in the Sobolev space $H^l(\R^d; N)$ for integer $l>d/2+1$ is established in \cite{McGahagan2007}.
One can see by scaling that $\Hdot{d/2}$ is critical for the problem, and significant work
	has been done on proving global well-posedness in this and other critical spaces in the special case 
	when the target is the sphere $N=\Sp^2 $
	\cite{Benjenaru2008,Benjenaru2009,Benjenaru2011,Ionescu2007}.

The case of the sphere is particularly attractive for two reasons. 
First, given the usual embedding $\Sp^2 \subset \R^3$, the \Sch{} maps equation becomes quite explicit.
In this framework, the complex structure at the point $u$ is simply given by the cross product in $\R^3$, $Jw = u\times w$.
The derivative term is calculated to be $\sum_k D_k \partial_k u = \Delta u + |\nabla u|^2 u$, where $\Delta$ and $\nabla$ are the Laplacian and 
	gradient operators for functions from $\R^d$ to $\R^3$.
The \Sch{} maps equation thus becomes,
\begin{equation}
	\label{eqn:1:heis}
	u_t = u \times (\Delta u + |\nabla u|^2u), \;\; x \in \R^d, \;\; u(x) \in \Sp^2 \subset \R^3.
\end{equation}

The second reason this case of the \Sch{} maps
equation is appealing is that it is physically relevant.
Equation \eqref{eqn:1:heis} is used to describe the dynamics of ferromagnetic spin systems, and is known in the physics community as the Heisenberg model.
It is a special case of the equation,
\begin{equation}
	\label{eqn:1:ll}
	u_t = (\alpha +  \beta u \times) (\Delta u + |\nabla u|^2u), \;\; x \in \R^d,
\end{equation}
which is the Landau-Lifshitz-Gilbert equation and is used to study the direction of magnetism in a solid.
(The survey article \cite{Lakshmanan2011} discusses the physical relevance of these equations.)
The equation \eqref{eqn:1:ll} corresponds
	precisely to the \gll{} equation \eqref{eqn:1:mf} in the case of maps $u: \R^d \rightarrow \Sp^2$.
The work on small data existence and uniqueness in a critical space for the \Sch{} maps
equation
	in this case of the sphere 
	culminated 
 	in \cite{Benjenaru2011}, which furnished a global critical small data well-posedness result in the Sobolev space $\Hdot{d/2}$.

A large body of work has been devoted to the sphere problem when the domain 
is $\R^2$.
The critical space is $\Hdot{1}$, so the problem in this dimension is energy critical.
It is also tractable to study because one can make an \emph{equivariant ansatz} and thereby study a sub-problem of the flow as a whole.
The equivariant ansatz involves studying solutions of the form
	the form 
	$u(r, \theta) = e^{m\theta R}f(r)$
	where $f(r) \in \R^3$, $m \in \Z$, and $R$ is the generator of rotations about the $z$-axis and given by the matrix
\[
	R= \begin{pmatrix}
		 0	&1	&	0\\
		 -1	&0	&	0\\
		 0	&0	&	0
	\end{pmatrix}.
\]
The overall picture that has developed can be described in terms of the harmonic maps,
	which have finite energy in this context,
	and whose existence is generally seen  as barrier to global well-posedness.
In the case of radial maps, $m=0$, there are no non-trivial harmonic maps and
	a global existence result for arbritrarily sized data in $H^2$ has been established  \cite{Gustafson2010}.
In the case when $m=1$, the lowest energy level of the non-trivial harmonic maps is $4\pi$;
	for initial data with energy strictly smaller than this, global existence has been shown to hold \cite{Benjenaru2013}.
On the other hand, in \cite{Merle2013}, a set of initial data with energy arbitrarily close to $4\pi$ is constructed which generates finite time blow up solutions.
	(This paper resolved the long standing question of whether finite energy initial data could lead to finite time blow up.)
Finite time blow up solutions are also constructed in \cite{Perelman2013}.
For $m\geq 3$, it has been shown that if the  initial data has energy close to that of the harmonic maps then the solution is, in fact, global \cite{Gustafson2010}.

Still in dimension 2, the equivariant ansatz can be made under the more general assumption that the target $N$ is a complex surface with an $\Sp^1$ symmetry.
This was originally done in \cite{Chang2000}, where a critical well-posedness theory for equivariant data small in $\Hdot{1}$ was developed.
Under the same equivariant ansatz, \cite{Germain2010} take a different approach than the Sobolev theory, and instead study
	the self-similar solutions of the flow.
These are solutions of the form $u(x,t) = \psi(x/\sqrt{|t|})$ for a profile $\psi$.
A family of such solutions with $C^\infty$ profiles is constructed, giving an example of 
regularity breakdown: these
	solutions are smooth at all times $t\neq 0$ but not smooth at $t=0$.
The study of these self similar solutions is supplemented with a global critical small data well-posedness theorem 
	in a Lorentz space that is shown to include the self-similar data.

When the dimensions of the range and domain are larger than two, but the same, it is still possible to formulate an equivariant ansatz,
	as will be shown in detail below.
For the case of the \Sch{} maps equation for maps $u:\C^n \rightarrow \CPn$, this equivariant ansatz is considered in \cite{Ding2008},
	where the existence of self-similar solutions is established.

The primary purpose of the present paper is to expand 
upon this previous work on the equivariant $\C^n$ to $\CPn$ case,
with a particular interest in establishing a global wellposedness theorem.
Our central result is a new equation for the dynamics in this case
    \eqref{eqn:1:sphere}.
This new equation is similar in structure to the \gll{} equation for maps
    to the sphere \eqref{eqn:1:ll}, and thus immediately opens up 
    the possibility
    of applying research ideas developed for the sphere problem to the
    present context.
Our global wellposedness result 
    in Section \ref{sec:wellposedness} in the case $n=2$
    is an example of this in practice.

\subsection{Overview of the results}

\subsubsection{The equivariant ansatz and derivation of the equation}

We consider maps $v: \C^n \rightarrow \CPn$, where $\CPn$ is equipped with the Fubini-Study metric, for $n\geq 2$.
The $n=1$ case is the usual problem of $\R^2$ to the sphere because $\CPn[1]$ with the Fubini-Study metric is isometric to $\Sp^2$.
In what follows $n$ is the complex dimension and $d=2n$ is the real dimension.

Recall that $\CPn$ can be viewed in terms of the homogeneous coordinates as
	points $(z_0,z_1,\ldots,z_n) \in \C^{n+1}$ under the identification
$
	[z_0,z_1,\ldots,z_n] = [\alpha z_0, \alpha z_1, \ldots, \alpha z_n],
$
for all $\alpha \in \C\backslash\{0\}$.
Given a complex isometry $A$ of $\C^n$ we can construct an isometry $\At$ of $\CPn$ by the formula
$
	\At [z_0,z_1,\ldots,z_n] = [z_0, A(z_1,\ldots,z_n)];
$
that is, we let $A$ act on the last $n$ coordinates in the homogeneous representation.
A map $v: \C^n \rightarrow \CPn$ is said to be \emph{equivariant} if $v( A z ) = \At v(z)$ for all isometries $A$ of $\Cn$ and all points $z\in \C^n$.
This ansatz is formally conserved by the flow.
This assumption is strong and, as we show, implies that $v$ is in fact of the form
$
	v(z) = 	v( (z_1,\ldots,z_n) ) = [z_0, f(r) z_1, \ldots ,f(r) z_n]
$
where $r = |z|$ and $f: \R^+ \rightarrow \C$.
We observe that for any $x \in \R$ we have
$
	v( (x,0,\ldots,0) ) = [z_0, f(r) x, 0 ,\ldots, 0] ;
$
or namely that,
\[
	v(\R^+ e_1) \subset \left\{ [z_0,z_1,0,\ldots,0] \;:\; z_0,z_1 \in \C \right\} \simeq \CPn[1],
\]
so the image of a real ray is contained in a complex line.
The Fubini-Study metric of $\CPn$ restricts to the Fubini-Study metric on this $\CPn[1]$,
	so in fact the image of $v(\R^+ e_1)$ is contained in a manifold isometric to $\Sp^2$.
The idea, now, is to parameterize this sphere in the usual 
embedding $\Sp^2 \subset \R^3$ and determine an equation on $u(r) = v(re_1) \in \Sp^2$.
From the equivariant ansatz we can recover $v$ from $u$.

By a computation we determine that the energy of $v$ is given in terms of $u:\C^n \rightarrow \Sp^2$ by the formula,
\begin{equation}
	\label{eqn:1:energy}
	\E(u) 
	= \frac{1}{2} \int_{\R^{2n}} \left(  | u_r |^2 + \frac{  u_1^2 + u_2^2 + (2n-2) | u - e_3 |^2  }{ r^2 }\right) dx,
\end{equation}
where $|u - e_3|$ is the Euclidean distance in $\R^3$ between $u$ and the north pole of the sphere $e_3$,
	and $|u_r|$ is the Euclidean norm in $\R^3$ of $u_r$. 
Observe that in the case $n=1$, we recover the usual energy for the equivariant $\R^2\rightarrow \Sp^2$ problem, as we would expect.
(See, for example, \cite{Benjenaru2013}, p.\ 2.)
For $n\geq 2$, one determines that any function $u$ with finite energy is continuous and has a limit as $r \rightarrow \infty$;
	by inspecting the energy one sees that this limit must be the north pole $e_3$.

The harmonic map heat flow, the \Sch{} maps equation, 
and the \gll{} equation for this equivariant case are now determined by calculating the variation of the energy.
We find that the \gll{} equation is given by,
\begin{equation}
	\label{eqn:1:sphere}
	u_t = (\alpha P + \beta u \times)\left( \der{^2 u }{r^2} + \frac{2n-1}{r} \der{u}{r}  + \frac{ 2n-2+u_3 }{r^2} e_3 \right),
\end{equation}
where 
 $P$ is the projection onto the tangent space $T_{u}\Sp^2$ and $u_3= \langle u,e_3 \rangle$.
The harmonic map heat flow corresponds to $\alpha=1$ and $\beta=0$; while the \Sch{}
 maps equation corresponds to $\alpha=0$ and $\beta=1$.
This representation of the problem appears to be new.
Its similarity to the corresponding equation  for maps to the sphere is precisely what makes it so useful:
	it immediately opens up the possibility of applying some of the techniques
        that have been developed for the case of the sphere to the present setting too.

By taking the stereographic projection from the north pole $f(r) =( u_1(r)+i u_2(r))/( 1+u_3(r) )$ we determine the stereographic representation of the problem,
\begin{equation}
	\label{eqn:1:conformal}
	f_t = (\alpha + \beta i) \left( f_{rr} - \frac{ 2 \of f_r^2 }{ 1+|f|^2} + \frac{2n-1}{r} f_r - \frac{2n-1}{r^2} f + \frac{1}{r^2} \frac{ 2 |f|^2 f }{ 1 + |f|^2  } \right),
\end{equation}
where the function here is a radial map $f: \R^{2n} \rightarrow \C$.
From this representation we see right away that the harmonic maps -- that is, the stationary solutions -- are given explicitly in this context by $f(r) = \alpha r$
for any $\alpha\in \C$.
In the terms of the sphere coordinates, the harmonic maps are given by a type of  stereographic projection,
\begin{equation}
	\label{eqn:1:harmonicsphere}
	u(r) =  \frac{1}{1+|\alpha|^2 r^2} \left( 2\Re (\alpha) r, 2\Im(\alpha) r , 1 -|\alpha|^2 r^2 \right).
\end{equation}
Again, this is consistent with the $n=1$ case, where the equivariant harmonic maps from $\R^2$ to $\Sp^2$ are known to be stereographic projections.
What is remarkable is that the analytic expressions for the harmonic maps are independent of $n$.
This seems to suggest that, from the perspective of the theory of harmonic maps, $\CPn$ is the natural higher dimensional analog of $\Sp^2$.

However there is a difference for $n\geq 2$: 
	observe that from \eqref{eqn:1:harmonicsphere} we have $\lim_{r \rightarrow \infty} u(r) = -e_3$, and so we find, by previous remarks on the energy,
	 that in this equivariant context all of the non-trivial harmonic maps have infinite energy.
	
\subsubsection{Self similar solutions}

After deriving the equation describing the dynamics, 
we first study the self similar solutions of the problem, which are given by $u(r,t) = \psi(r / \sqrt{t})$
	for a profile $\psi(r) = u(r,1)$.
By substituting this ansatz into \eqref{eqn:1:sphere} we determine the following ODE system on $\psi$:
\begin{align}
	0 &= \displaystyle (\alpha P + \beta u \times)\left( \frac{ \partial^2 \psi }{ \partial r } + \left( \frac{2n-1}{r} + \frac{r}{2} \right) \der{\psi}{r} + \frac{ 2n-2+\psi_3 }{r^2} e_3 \right),	
	\label{eqn:1:selfsim}
	\\
	\psi(0) & = e_3,\nonumber	\\
	\psi'(0) &= v = (v_1,v_2,0) \in T_{e_3} \Sp^2. \nonumber
\end{align}
As mentioned previously, the self similar solutions for the \Sch{} 
maps equation in this equivariant setting have already been studied in \cite{Ding2008}.
However using the representation \eqref{eqn:1:selfsim} we are able to
simplify the analysis significantly.
We are also able to extend the analysis by gaining more information
on the convergence of the self-similar profile, and by
treating the general \gll{} equation as well as the \Sch{} maps equation.

\begin{thm}
	\label{thm:1}
	Fix $\alpha\geq 0$ and $\beta \in \R$.
	For every $v \in T_{e_3} \Sp^2$ there is a unique global solution to \eqref{eqn:1:selfsim}.
	The solution is smooth for $r>0$. 
	In the non-trivial case, when $v\neq 0$, the solution has the following properties:
	\begin{enumerate}
		\item	For all $r>0$, $\psi(r) \neq e_3$.
		\item	If $\alpha>0$ then $|\psi_r| \lesssim 1/r^3$. If $\alpha=0$ then $r |\psi_r| \rightarrow 0$ as $r \rightarrow \infty$.
		\item	If $v \neq 0$, there exists a point $\psi_\infty  \in \Sp^2$, $\psi_\infty \neq e_3$, such that $\lim_{r \rightarrow \infty} \psi(r) = \psi_\infty$.
			Consequently, $\E(\psi)=\infty$.
		\item	The limit $\psi_\infty$ depends continuously on $v$;	
				in particular, $\lim_{v \rightarrow 0} \psi_\infty = e_3$.
	\end{enumerate}
\end{thm}

Because of the convergence, we see that $u(r,t) = \psi(r/\sqrt{t})$ is a solution of the \gll{} flow
	corresponding to the initial data $u(r,0) \equiv \psi_\infty$.

Notice that in the case $\alpha>0$ -- that is, when there is some dissipation
-- we are able to prove
	faster convergence to 0 of $\psi_r$.
In the case of the \Sch{} maps equation ($\alpha=0$)
the rate of convergence of $\psi_r$ is insufficient to guarantee the convergence of $\psi$,
	so an additional argument is needed.

\subsubsection{Global critical wellposedness}

We next
illustrate how methods for proving wellposedness of the 
    \Sch{} maps equation for the sphere may be adapted to prove
    wellposedness of \eqref{eqn:1:sphere}.
We specifically adapt the Hasimoto transformation method from \cite{Chang2000}.
For a smooth solution $u(r,t)$ of \eqref{eqn:1:sphere} and a fixed time $t$,
    the map $r \mapsto u(r,t)$ defines a curve on $\Sp^2$ starting at $e_3$.
Choose any element $e \in T_{e_3} \Sp^2$ and consider the parallel
transport $e(r)$ of this curve along $r \mapsto u(r,t)$.
Because the tangent space at the point $u(r,t)$ of the sphere is 
    two-dimensional, it is spanned by $e(r)$ and $Je(r) = u\times e(r)$.
We may therefore define a complex valued function $q$ by the formula,
\begin{equation}
    \Re(q) e(r) + \Im(q) J e(r) = q e(r) = u_r,
    \label{eqn:1:hasimoto}
\end{equation}
precisely as in \cite{Chang2000}.
This equation is known as the \emph{Hasimoto transformation}.
It is chosen so that the function $q$ will satisfy a `nice' nonlinear
\Sch{} equation; namely, an equation where the non-linearity does
not contain derivatives.
We derive the equation on $q$ for all $n$, and in the case
    $n=2$ -- that is, for the equivariant \gll{} maps equation 
from $\C^2$ to $\CPn[2]$ -- we provide the necessary estimates
to prove the following small-data critical global wellposedness result.

\begin{thm}
	\label{thm:2}
	Fix $p \in [1,2]$.
        Define $r$ by $1/r = 1/2 - 1/6p$
	and the spaces $X$ and $X_0$ by the norms
	\[
		\| q \|_X = \| \nabla q \|_{L^{3p}_t L^{r}_x}
		\;\;\text{ and }\;\;
		\| q \|_{X_0} = \| e^{it\Delta} (aq) \|_{X},
	\]
	where $a(x) = x_1/r$.
        There exists $\epsilon>0$ such that if 
            $u_0: \R^{2n} \rightarrow \Sp^2$ is radial,
            $q_0$ is defined by \eqref{eqn:1:hasimoto}, and
            $\| q_0 \|_{X_0} \leq \epsilon$,
        there is unique global solution of the 
        \gll{} equation \eqref{eqn:1:sphere} for $\beta>0$
            for $n=2$ with the derivate term $q$ in the space $X$.
\end{thm}

Some remarks.
\begin{itemize}
\item
This is, to the best of our knowledge, the first global wellposedness
result for the \Sch{} maps equation where the target manifold has complex
dimension greater than one.
\item
The space $X$ is at the scaling level of the equation.

\item
Because $(3p,r)$ is an admissible exponent pair for the Strichartz
estimates for the \Sch{} equation, we have
$
    \| q \|_{X_0} \lesssim \| \nabla (aq) \|_{L^2} \lesssim \| \nabla q \|_{L^2}
    $
    and hence data $q_0$ whose derivative is small in $L^2$
    are included in the wellposedness result.

\item
For $n>2$ we are unable to provide the estimates to close the argument 
in an elementary way.
A global wellposedness result for arbitrary $n$,
proved using the Hasimoto transform or another method adapted from the 
research on the \Sch{} maps equation for the sphere, would be
very satisfactory.
\end{itemize}

\subsubsection{The `real' heat flow case}

We finally study
an interesting sub-problem of the general equation
\eqref{eqn:1:sphere} corresponding to
the harmonic map heat flow with an additional
condition on the inital data.
Recall that for the linear heat equation, if one starts with real valued data then the solution will be real valued for all time.
On the other hand, if one starts the linear \Sch{} equation with real valued data then the solution will, in general, be complex valued for future times.
This shows that in the heat flow case there is a lower dimensional sub-problem when one restricts to real valued data.

In our context, the analogous fact
is that if one starts the harmonic map heat flow \eqref{eqn:1:sphere}
with initial data valued in a great circle passing through the north pole,
	the solution will continue to be valued on the same great circle for future times. 
For the \gll{} flow this is not true: the solution will spread out to the whole sphere.
For the harmonic map heat flow one can thus fix a great circle and consider the problem for initial data valued on that circle.
One expects the analysis of this sub-problem to be easier as the dimension of the problem is reduced.
However, because both the harmonic maps and the self similar solutions are solutions of this type,
	 it is still an interesting case to consider.

By paramaterizing the great circle by its spherical distance from the north pole, one finds that the `real' heat flow is given by the PDE,
\begin{equation}
	\label{eqn:1:realheatflow}
	g_t = \frac{ \partial^2 g }{ \partial r^2 } + \frac{2n-1}{r} \der{g}{r}  + \frac{ \eta(g)  }{ r^2 },
\end{equation}
where $\eta(g) = \sin(2g) + (2n-2) \sin(g)$.
Equations of this type, which arise in the study of the equivariant harmonic map heat flow 
	on spherically symmetric manifolds,  have been extensively studied \cite{Germain2011,Germain2017}.
There is a general theorem which, based on the structure of $\eta$, classifies the PDE
	into a uniqueness regime or a non-uniqueness regime.
Our primary purpose here is to show that for $n=2$ -- that is, the problem of maps from $\C^2$ to $\CPn[2]$ --  the PDE \eqref{eqn:1:realheatflow} is a \emph{borderline case}
	for this classification theorem.
We find that the dynamics of the PDE share some of features of the uniqueness regime, and some of the features of the non-uniqueness regime, but ultimately
	that non-uniqueness holds.

\begin{thm}
	
	(i) For $n=2$ there is a weak non-constant solution of \eqref{eqn:1:realheatflow} corresponding to the initial data $g_0(r) \equiv \pi$.
		This solution is distinct from the constant solution $g(r,t) \equiv \pi$.

	(ii) In the case $n\geq3$, for each initial data in $L^\infty$ and each $T>0$,
			there is at most one solution of \eqref{eqn:1:realheatflow} in $L^\infty([0,T],L^\infty)$.
\end{thm}

\section{The equivariant ansatz and derivation of the equation}

\label{sec:derivation}

\subsection{The equivariant ansatz}

We consider maps $v: \C^n \rightarrow \CPn$.
In order to rigorously describe the equivariant ansatz, we recall more carefully the construction of $\CPn$.
One begins with vectors $z=(z_0,z_1,\ldots,z_n)\in\C^{n+1} \backslash \{0\}$ and first identifies points $z \sim \lambda z$
where $\lambda \in \R \backslash \{0\}$.
The resulting equivalence classes can be identified 
	with points on the sphere $\Sp^{2n+1} \subset \C^{n+1}$.
This sphere has the usual metric induced from $\C^{n+1}$.
Now one defines the equivalence relation $z \sim e^{i\theta} z$ for $\theta \in \R$, and defines $\CPn = \Sp^{2n+1}/\sim$.
The Fubini-Study metric is the metric induced from $\Sp^{2n+1}$.

To make the equivariant ansatz, we first construct a special class of isometries on $\CPn$ in the following way.
Take any complex isometry $A$ of $\C^n$, and define $\hat{A}:\Cn[n+1] \rightarrow \Cn[n+1]$ by,
\[
	\hat{A}( z_0,z_1,\ldots,z_n) = (z_0, A(z_1,\ldots,z_n));	
\]
that is, $A$ acts on the last $n$ coordinates of a point in $\Cn[n+1]$.
If $A$ is a complex isometry of $\Cn$, then $\hat{A}$ is clearly a complex isometry of $\Cn[n+1]$.
Now define a map $\tilde{A}$ on $\CPn$ through the homogeneous coordinates by,
\begin{equation}
	\label{eqn:2:iso}
	\tilde{A}[z_0,z_1,\ldots,z_n] 
		= [\hat{A}(z_0,z_1,\ldots, z_n)]
		= [z_0, A(z_1,\ldots, z_n)].
\end{equation}
The map $\tilde{A}$ is well defined because $A$ commutes with complex scalar multiplication.

\begin{lemma}
        If $A$ is a complex isometry of $\Cn$ then 
        $\tilde{A}$ defined by \eqref{eqn:2:iso} is an
		isometry of $\CPn$.
\end{lemma}

\begin{proof}

	We have
	\begin{align*}
		d_{\CPn}( \tilde{A}[v], \tilde{A}[w] )
                    &= d_{\CPn}( [\hat{A}v], [\hat{A}w] )						        = \min_{\alpha,\beta \in [0,2\pi] } d_{\Sp^{2n+1}}
                        ( e^{i\alpha}\hat{A}v, e^{i\beta }\hat{A}w )	\\
		    &= \min_{\alpha,\beta \in [0,2\pi] } 2 
                        \arcsin\left( \frac{1}{2} d_{\C^{n+1}} 
                        (e^{i\alpha}\hat{A}v,e^{i\beta}\hat{A}w)  \right) \\
		    &= \min_{\alpha,\beta \in [0,2\pi] } 2 
                        \arcsin\left( \frac{1}{2} d_{\C^{n+1}} 
                        (e^{i\alpha}v,e^{i\beta}w)  \right) 
                        =  d_{\CPn}( [v],[w] ),
	\end{align*}
	where in the second to last equality we used that $\hat{A}$ 
        commutes with $e^{i\theta}$ and that $\hat{A}$ is an isometry of $\C^{n+1}$.
\end{proof}

We say a map $v: \Cn \rightarrow \CPn$ is \emph{equivariant} if $v(Az) = \tilde{A} v(z)$ for all complex isometries $A$ of $\Cn$.
We now show that this assumption implies a strong rigidity on $v$.
Take any $z\in\Cn$ and write $v(z) = [w_0,w]$ for some $w_0\in\C$ and $w\in\C^n$.
Now consider any isometry $A$ that fixes $z$.
By the equivariant ansatz and $Az=z$ we have,
\[
	[w_0, Aw] = \tilde{A} u(z) = u(Az) = u(z) = [w_0,w],
\]
which implies that $Aw = w$, so $A$ also fixes $w$.
Because $A$ is an arbitrary isometry that fixes $z$, we must in fact have $w = f(z) z $ for some $f(z) \in \C$,
	and hence $v(z) = [w_0,f(z)z]$ for all $z$. Moreover, we have,
\[
	[w_0,f(Az) Az ] = v(Az) = \tilde{A}v(z) = [w_0,A( f(z) z) ] = [w_0, f(z) Az ],
\]
so $f(Az)=f(z)$. Because this holds for all isometries $A$, $f(z)$ is in fact a radial function
and hence,
\begin{equation}
\label{eqn:2:nueqn}
	v(z) = [w_0,f(|z|)z],
\end{equation}
for some function $f: \R^+ \rightarrow \C$.

We now observe that if $r\in\R^+$ then $v(r e_1)=[w_0,f(r)r,0,\ldots,0]$.
In other words
\[
	v(\R^+ e_1) \subset \{ [w_0,w_1,0,\ldots,0 ] \;:\; w_0,w_1 \in \C \} \simeq \CPn[1].
\]
The Fubini-Study metric  on $\CPn$ restricts to the Fubini-Study metric on $\CPn[1]$,
	and so this $\CPn[1]$ is isometric to the sphere $\Sp^2$.
Moreover, the complex structure of $\CPn$ restricts to the standard complex structure of $\CPn[1]$.
In the usual embedding $\Sp^2 \subset \R^3$ this is given, as is well known, by $J w = u \times w$ at the point $u\in \Sp^2$ and for all $w \in T_u\Sp^2$.
We next parameterize this sphere and determine an equation for the function $r \mapsto v(re_1) \in \Sp^2$.

\subsection{Derivation of the energy}

The isometric identification between $\CPn[1]$ (with the Fubini-Study metric) and $\Sp^2\subset \R^3$ (with the metric from the standard embedding)
	can be made through the isometric invertible map,
\begin{equation}
	\label{eqn:2:spcpn}
	\Sp^2 \ni (a_1,a_2,a_3) \mapsto \frac{1}{ \sqrt{2}  (1+a_3)^{1/2} } [1+a_3, a_1 + i a_2] \in \CPn[1],
\end{equation}
where in this case the north pole $e_3=(0,0,1)$ is mapped to the point 
$[1,0] \in \CPn[1]$.
In this identification the complex structure on $\CPn[1]$ is mapped to the standard complex structure on the sphere.
Given an equivariant map $v: \C^n \rightarrow \CPn$, we wish to write it in a form so that
	$v(re_1) \in \CPn[1]$
        has the representation $[1+a_3,a_1+ia_2,0,\ldots,0]$.
In fact, we can write $v$ in the form,
\begin{equation}
	\label{eqn:2:vparam}
	v(z) = \frac{1}{ \sqrt{2} (1+u_3)^{1/2} } \left[ 1+u_3, ( u_1+iu_2) \frac{z}{r} \right]
\end{equation}
for $u(r) = (u_1(r),u_2(r),u_3(r))$ satisfying $|u|_{\R^3}=1$.
When we substitute $z=re_1$ we recover essentially the representation in \eqref{eqn:2:spcpn},
	and hence $u$ parameterises the sphere in the correct, isometric, way.

(To see that $v(z) = [w_0,f(r)z]$ in \eqref{eqn:2:nueqn} can be written as in \eqref{eqn:2:vparam}, 
	observe that by scaling we can assume that $(w_0,g(r)z)\in\Sp^{2n+1}$,
	which means $|w_0|^2 + |g(r)|^2r^2 = 1$.
We can also assume by scaling that $w_0>0$.
This means, in fact, that $w_0 \in [0,1]$, and 
	hence there is a unique $u_3(r) \in [-1,1]$ such that $\sqrt{2}(1+u_3(r))^{1/2} = w_0$.
We then define $u_1 + i u_2 = r g(r) \sqrt{2}(1-u_3)^{1/2}$, and substuting this in gives the representation above.
The condition $|w_0|^2 + |g(r)r|^2=1$ translates into $|u|_{\R^3} = 1$.)

\begin{prop}
	The energy is given in the $u$ coordinates by,
	\begin{equation}
		\label{eqn:2:energy}
		\E(v) = \frac{1}{2}  \int_{\R^{2n}} |dv|^2 dx =
		\frac{1}{2} \int_{\R^{2n}} \left[ |u_r|^2 +  \frac{1}{r^2} \left[ 1-u_3^2  + 2(2n-2)(1-u_3) \right]  \right] dx.
	\end{equation}
\end{prop}

\begin{proof}
In order to calculate the energy density $|dv|^2$ of $v(z)$ we have to fix a basis for $T_z\C^n$,
	which will be $2n$ dimensional,
	and calculate first derivatives of $v$ with respect to this basis. 
For concreteness we view $v$ as being valued in the sphere $\Sp^{2n+1}$,
\begin{equation}
	\label{eqn:2:vinsphere}
	v(z) = \frac{1}{ \sqrt{2} (1+u_3)^{1/2} } \left( 1+u_3, ( u_1+iu_2) \frac{z}{r} \right) \in \Sp^{2n+1} \subset \C^{n+1},
\end{equation}
and perform the computation there.
The only adjustment needed to be made is as follows.
Given a point $p \in \Sp^{2n+1}$, all points $e^{i\theta}p$ are mapped to the same point $[p]\in \CPn$.
By differentiating with respect to $\theta$, it is apparent that in $T_p\Sp^{2n+1}$ the tangent direction $ip \in T_p\Sp^{2n+1}$ is contracted
	under the identification $p \sim e^{i\theta }p$.
Hence when calculating derivatives at the level of $\Sp^{2n+1}$ we take usual Euclidean derivatives in $\C^{n+1}$,
	project onto $T_p\Sp^{2n+1}$, and then factor out the real subspace spanned by $ip$.
In fact, that the last two parts of this process amount to taking the complex projection,
\begin{equation}
    \label{eqn:2:projection}
	P w = w - \langle w,p\rangle_{\C^{n+1}} p
\end{equation}
of derivative terms $w$.
We have, of course, $|Pv|^2 = |v|^2 - | \langle v, p \rangle |^2$.

Let $\partial/\partial z_k$ and $\partial/\partial \zkb$ be the usual basis for $T_z \C^n$.
For any vector $w_0 \in \Cn$ define
	 $\partial/\partial w_0 = \sum_{m=1}^n w_0^m \partial/\partial z_m$ and
	 $\partial/\partial \overline{w_0} = \sum_{m=1}^n \overline{w_0^m} \partial/\partial \overline{z_m}$.
If $\{ w_k \}_{k=1}^n$ is an orthonormal basis of $\Cn$ then the derivatives $\{ \partial/\partial w_k, \partial/\partial \overline{w_k} \}$ are an orthogonal basis for the tangent space and so,
	by the expression for $|dv|^2$ local in coordinates,
\begin{equation}
	\label{eqn:2:dvcoords}
	|dv|^2 = 4 \sum_{k=1}^n 
		\left| P \frac{\partial v}{\partial w_k} \right|^2 + 
		\left| P \frac{\partial v}{\partial \overline{w_k}} \right|^2.
\end{equation}
One verifies the formulas at the point $z\in \Cn$,
\begin{align}
	\frac{\partial r }{\partial w_0 } 		&= \frac{ \langle w_0, z \rangle}{ 2r};		&
	\frac{\partial r }{\partial \overline{w} } 	&= \frac{ \langle z, w_0 \rangle}{2r};	&
	\frac{\partial z }{\partial w_0 }		&= \frac{w_0}{|w_0|};				&
	\frac{\partial z }{\partial \overline{w_0} }	& = 0.	
\end{align}
We then set $w_1=z/|z|$ and define $w_2(z),\ldots,w_n(z)$ locally so that that $\{w_k(z)\}_{k=1}^n$ is 
	an orthonormal basis of $\C^n$ for each $z$.
In this setup, $w_1$ is the radial direction and $w_k$ derivatives for $k\geq 2$ will be independent of radial terms.

Hence for $k \geq 2$ we compute and find,
\[
	\frac{ \partial v }{ \partial w_k } = \frac{1}{ \sqrt{2} (1+u_3)^{1/2} } \left( 0, (u_1+iu_2) \frac{w_k}{r} \right)	\;\;\text{ and }\;\;
	\frac{ \partial v }{ \partial \overline{w_k} } = 0.
\]
We see from \eqref{eqn:2:vinsphere} that $\partial v/\partial w_k$ is complex orthogonal to $v$ and so,
\[
	\left| P \frac{ \partial v }{ \partial w_k } \right|^2 = 
	\left|  \frac{ \partial v }{ \partial w_k } \right|^2 = 
		\frac{1}{ 2 (1+u_3) } \frac{ u_1^2 + u_2^2 }{ r^2 } =
		\frac{1-u_3}{ 2r^2 },
\]
where in the step we used $u_1^2 + u_2^2 + u_3^2=1$.

We now differentiate with respect to $w_1$ and $\overline{w_1}$.
In this case the radial terms will also be differentiated.
We note, however, that when differentiating that we can ignore the scaling term $1/(\sqrt{2}(1+u_3)^{1/2})$: when this is differentiated
	we simply get a scalar multiple of $v(z)$, which disappears under the projection
        \eqref{eqn:2:projection}.
Hence,
\begin{align*}
	P \frac{ \partial v }{ \partial w^1 }	
		&= \frac{1}{ \sqrt{2} (1+u_3)^{1/2}} P \frac{\partial }{ \partial w^1 }
			\left[ 1+u_3, (u_1+i u_2) \frac{z}{r} \right]		\\
		&= \frac{1}{ \sqrt{2} (1+u_3)^{1/2}} P 
			\left( u_3' \frac{1}{2}, (u_1'+i u_2') \frac{z}{2r} - (u_1+iu_2) \frac{ z }{ 2r^2 } \right), 
\end{align*}
and similarly,
\begin{align*}
	P \frac{ \partial v }{ \partial \overline{w^1} }	
		&= \frac{1}{ \sqrt{2} (1+u_3)^{1/2}} P 
			\left( u_3' \frac{1}{2}, (u_1'+i u_2') \frac{z}{2r} + (u_1+iu_2) \frac{ z }{ 2r^2 } \right). 
\end{align*}
The difference in sign gives rise to the simplification
\begin{align*}
	\left| 	P \frac{ \partial v }{ \partial w^1 } \right|^2 
	+ \left| P \frac{ \partial v }{ \partial \overline{w^1} } \right|^2
		&= \frac{1}{ 4 (1+u_3) } \left[  \left| 
				P \left( u_3' , (u_1'+i u_2') \frac{z}{r} \right) \right|^2 +
				\left| P \left( 0, (u_1+iu_2) \frac{ z }{ r^2 } \right)  \right|^2
			 \right].	
\end{align*}
Finally, a computation using the relations $u_1^2+u_2^2+u_3^2=1$ and $u_1u_1'+u_2u_2'+u_3u_3'=0$ reveals that
\begin{align*}
	\left| P \left( u_3' , (u_1'+i u_2') \frac{z}{r} \right) \right|^2 	
		&= \left|  \left( u_3' , (u_1'+i u_2') \frac{z}{r} \right) \right|^2 	\\
		&\hspace{1cm}
			- \left| \left\langle 
				\left( u_3' , (u_1'+i u_2') \frac{z}{r} \right),
				\frac{1}{ \sqrt{2} (1+u_3)^{1/2} } \left( 1 + u_3, (u_1+iu_2) \frac{z}{r} \right) 
			\right\rangle \right|^2		\\
		&= |u_r|^2 - \frac{1}{2(1+u_3)}  \left| (1+u_3) u_3' + (u_1-i u_2)(u_1'+i u_2') \right|^2	\\
		&= |u_r|^2 - \frac{1}{2(1+u_3)}  \left|u_3' + i (u_1 u_2' - u_1' u_2) \right|^2	\\
		&= |u_r|^2 - \frac{1}{2(1+u_3)}  \left[ (u_3')^2 + u_1^2 ( u_2')^2 + (u_1')^2 u_2^2 -2  u_1u_1'u_2u_2' \right]	\\
		&= |u_r|^2 - \frac{1}{2(1+u_3)}  \left[ (1- u_3^2) | u_r |^2 \right]	
		= \frac{1+u_3}{2} |u_r|^2,
\end{align*}
and,
\begin{align*}
	\left| P \left( 0, (u_1+iu_2) \frac{ z }{ r^2 } \right)  \right|^2	
		&=	\left| P \left( \frac{-1-u_3}{r}, 0 \right)  \right|^2 
                \\&
                =	
                \frac{1}{r^2} \left[ (1+u_3)^2 - \frac{1}{ 2(1+u_3) } (1+u_3)^4 \right]	
		=	\frac{(1+u_3)(1-u_3^2)}{ 2r^2}.
\end{align*}
We have, then, by substituting these expressions into \eqref{eqn:2:dvcoords},
\[
	|du|^2 = \frac{|u_r|^2}{2} + \frac{1}{r^2} \left[ \frac{ 1-u_3^2 }{2} + 2(n-1)(1-u_3) \right],
\]
and then 
\[
	\E(v) = \int_{\R^{2n}} |dv|^2 dx = 
	\frac{1}{2} \int_{\R^{2n}} \left[ |u_r|^2 +  \frac{1}{r^2} \left[ 1-u_3^2  + 2(2n-2)(1-u_3) \right]  \right] dx,
\]
which completes the computation.
\end{proof}

By the relations $1-u_3^2=u_1^2 + u_2^2$ and 
$
	|u-e_3|^2 = u_1^2 + u_2^2 + (u_3-1)^2 = 2(1-u_3),
$
we can equivalently write the energy as in an $L^2$ form as,
\begin{equation}
	\E(v) = \frac{1}{2} \int_{\R^{2n}} \left[ |u_r|^2 +  \frac{1}{r^2} \left[ u_1^2 + u_2^2  + (2n-2) |u-e_3|^2 \right]  \right] dx.
	\label{eqn:2:energyL2}
\end{equation}
With this representation we determine the following result.
\begin{prop}
	There holds $\| u_r \|_{L^2}^2 \lesssim \E(u) \lesssim \|u_r\|_{L^2}^2$.
\end{prop}
\begin{proof}
	The lower bound is obvious.
	For the upper bound, we observe that 
        $u_1^2+u_2^2 \leq | u - e_3 |^2 $ and hence that,
	\[
		\E(u) \leq \frac{1}{2} \left(  \| u_r \|_{L^2}^2 + (2n-1) \left\| \frac{ u-e_3 }{ r } \right\|_{L^2}^2 \right),
	\]
	and the result follows from the Hardy inequality
        $\| \phi /r \|_{L^2} \lesssim \| \phi_r \|_{L^2}$ for 
        functions $\phi: \R^d \rightarrow \R^3$
        (see Theorem 
    \ref{thm:app:hardy1} in the Appendix).
\end{proof}

\subsection{Variation of the energy, and the flow PDEs}

In order to find the PDEs corresponding to the harmonic map heat flow, the \Sch{} 
maps equation, and the \gll{} equation, we need to calculate the variation of the energy,
	given by the formula
\[
	\int_{\R^{2n}} \langle \tau(u),w \rangle_{T_u\Sp^2} dx = -
		 \left. \frac{d}{d\epsilon} \right|_{\epsilon=0} \E( u + \epsilon w ),
\]
	for all radial maps $w: \R^{2n} \rightarrow  \in T \Sp^2$ such that $w(r) \in T_u\Sp^2$.

\begin{prop}
We have,
\[
	\tau(u) = P_{T_u\Sp^2} \left( \der{^2 u}{r^2} + \frac{2n-1}{r} \der{u}{r} + \frac{ 2n-2+u_3 }{r^2} e_3  \right).
\]
\end{prop}

\begin{proof}
Using the representation \eqref{eqn:2:energy} we find for $w \in T_u \Sp^2$,
\begin{align*}
	\left. \frac{d}{d\epsilon} \right|_{\epsilon=0} \E( u + \epsilon w )
		&= \frac{1}{2} \int_{\R^{2n}} 2  \langle u_r, w_r \rangle_{\R^3}   + \frac{1}{r^2} \left[   - 2  u_3 w_3 + 2(2n-2)( -  w_3 ) ) \right]	\\
		&= - \int_{\R^{2n}} \left\langle u_{rr} + \frac{2n-1}{r}u_r, w \right\rangle_{\R^3} + \frac{1}{r^2}  \left\langle (2n-2+u_3)e_3, w \right\rangle_{\R^3} dx	\\
		&= - \int_{\R^{2n}} \left\langle
				P_{T_u\Sp^2} \left( \der{^2 u}{r^2} + \frac{2n-1}{r} \der{u}{r} + \frac{ 2n-2+u_3 }{r^2} e_3 \right) ,w \right\rangle_{T_u\Sp^2} dx,
\end{align*}
and the formula follows.
\end{proof}

In general the harmonic map heat flow is given by $u_t = \tau(u)$, the \Sch{} 
maps equation is given by $u_t = J \tau(u)$, 
	where $J$ is the complex structure on the target, and the \gll{} equation
        is given by $u_t = (\alpha + \beta J) \tau(u)$ for $\alpha \geq 0 $ and $\beta \in \R$.
By the previous proposition, $\tau(u)$ is determined, while as discussed above, the complex structure in the $u$ coordinates is
	precisely the usual complex structure on the sphere. 
We are therefore ready to write down the flow PDEs.

\begin{defn}
	The  equivariant generalized Landau-Lifshitz (GLL) problem from $\C^n$ to $\CPn$ is the Cauchy problem
	for $u: \R^{2n} \rightarrow \Sp^2$ given by 
	\begin{align}
		u_t(r,t) &= \displaystyle (\alpha P + \beta u\times) \left( \der{^2 u}{r^2} + \frac{2n-1}{r} \der{u}{r} + \frac{ 2n-2+u_3 }{r^2} e_3  \right), 	
		\label{eqn:2:sphere}\\
		u(r,0) &= u_0(r),\text{ with }u_0(0) = 0, \nonumber
	\end{align}
	for $\alpha \geq 0 $ and $\beta \in \R$.
	The case  $\alpha = 1$ and $\beta=0$ is the harmonic map heat flow.
	The case $\alpha=0$ and $\beta=1$ is the \Sch{} maps equation.
\end{defn}

Note by re-scaling time we can always assume that $\alpha^2+\beta^2=1$, 
which we do from now on.

By taking the stereographic projection $f(r) = (u_1 + iu_2)/(1+u_3)$,
with inverse given by
\begin{equation}
	\label{eqn:2:stereo}
	(u_1,u_2,u_3) = \frac{1}{1+|f|^2} ( 2 \,\Re f, 2\, \Im f, 1-|f|^2 ),
\end{equation}
we can determine the \emph{stereographic representation} of the problem.
With this stereographic projection, the north pole is mapped to the origin.

\begin{prop}
	The \gll{} equation is given in the stereographic coordinates by
	\begin{equation}
		\label{eqn:2:conformal}
		f_t = (\alpha + i \beta) \left[ f_{rr} - \frac{ 2\of f_r^2 }{ 1+|f|^2 } + \frac{2n-1}{r} f_r - \frac{2n-1}{r^2} f + \frac{1}{r^2} \frac{ 2 |f|^2 f }{ 1+|f|^2 } \right].
	\end{equation}
\end{prop}

The proof involves substituting the expression for the stereographic projection \eqref{eqn:2:stereo} into the PDE \eqref{eqn:2:sphere} and computing;
	we omit this standard computation.

\subsection{Classification of the harmonic maps in this context}

The equivariant harmonic maps from $\Cn$ to $\CPn$ are the time independent solutions of \eqref{eqn:2:sphere}.
Because the PDE has one space dimension, the time independent problem is an ODE.
In all, $\phi$ is harmonic if and only if,
\begin{equation}
	\label{eqn:2:harmapode}
	0 = \phi \times \left( \frac{ d^2 \phi }{ d r^2} + \frac{2n-1}{r} \frac{ d \phi }{ d r }
			+ \frac{ (2n-2) + \phi_3  }{ r^2 } e_3 \right),
\end{equation}
with the boundary conditions given by $\phi(0)=e_3$ and $\phi'(0)=v = (v_1,v_2,0) \in T_{e_3} \Sp^2 $.
Writing the harmonic function $\phi$ in the stereographic coordinates as $g$, the ODE is,
\[
	0 = g_{rr} - \frac{ 2 \overline{g} g_r^2 }{ 1+|g|^2 } + \frac{2n-1}{r^2}g_r - \frac{2n-1}{r^2} g + \frac{1}{r^2} \frac{2 |g|^2 g }{1+|g|^2}
\]
and the boundary conditions are $g(0)=0$ and $g_r(0) = v_1 + i v_2$.
Remarkably, we can solve this ODE explicitly with the linear function $g(r) = (v_1+iv_2)r$.
Moreover, because it is an ODE for which we have a uniqueness theory, $g(r) =(v_1+iv_2) r$ is the unique solution.
(See the Theorem \ref{thm:app:wellposedness} in the
Appendix for a local well-posedness theory for ODE of this type.)
Using the stereographic projection we can write the harmonic map in the sphere coordinates as
\begin{equation}
	\label{eqn:2:harmap}
	\phi(r) = \frac{1}{ 1 + |v|^2 r^2 } \left( 2 r v_1, 2 r v_2, 1-|v|^2 r^2  \right) =  \frac{1}{ 1 + |v|^2 r^2 } \left( 2r v + (1-|v|^2 r^2)e_3 \right);
\end{equation}
in fact, $\phi$ is just a version of the stereographic projection itself.
This is consistent with the well-known fact that the harmonic maps in the sphere ($n=1$) case are stereographic projections;
	what is interesting is that when $n$ is incremented in the ODE \eqref{eqn:2:harmapode}, the new terms still cancel under this expression.

Qualitatively speaking, the harmonic maps in our context are quite simple: they start, when $r=0$, at the north pole and, as $r$ increases, move monotonically away from the
	north pole, converging to the south pole in the limit $r \rightarrow \infty$.
By way of comparison, in the case of equivariant harmonic maps from the $d$-dimensional ball $B^d$ to $\Sp^d$ the situation is different \cite{Jager1983}.
For $3 \leq d \leq 6$
	the harmonic maps oscillate about the south pole, while for $d \geq 7$ the harmonic maps approach the south pole monotonically, as here. 
In general one finds that the equivariant harmonic maps usually fall into either an oscillatory regime or a monotonic regime \cite{Germain2010}.

Finally, we note that while the expressions above for the harmonic maps are independent of $n$, there is a difference when $n \geq 2$.
In the case of the sphere, $n=1$, the energy of the stereographic projection  is $4\pi$.
(This may be verified by substituting \eqref{eqn:2:harmap} into the energy \eqref{eqn:2:energy} with $n=1$, or by consulting \cite{Benjenaru2013}.)
However, for $n\geq2$ the energy is infinite.
To see this it is sufficient to observe that
$
	\lim_{r \rightarrow \infty} \phi(r) = -e_3
$
and to use the following Lemma.

\begin{lemma}
	\label{lemma:2:infenergy}
	Suppose that $\E(u)<\infty$ and $n\geq 2$. Then $\lim_{r \rightarrow \infty} u(r)$ exists and equals $e_3$.
\end{lemma}
\begin{proof}
	For any $r_2>r_1>0$ we have
	\begin{align*}
		|u(r_2)-u(r_1)| &= \left| \int_{r_1}^{r_2} u_r(r)dr \right| \leq \left( \int_{r_1}^{r_2} |u_r|^2 r^{2n-1} dr \right)^{1/2} 
			\left( \int_{r_1}^{r_2} \frac{1}{r^{2n-1}} dr \right)^{1/2}	
			\leq C \E(u) r_1^{-n+1},
	\end{align*}
	which, because $n\geq2$, shows that $\lim_{r \rightarrow \infty} u(r)$ exists.
	This means that in the energy \eqref{eqn:2:energyL2}, the right most term in the integrand,
	$
		(1/r^2) (2n-2)|u(r)-e_3|^2 r^{2n-1}
	$ 
	converges as $r \rightarrow \infty$.
	For the energy  to be finite, the limit must be 0.
	As $n\geq 2$, this implies that $\lim_{r\rightarrow\infty}u(r)=e_3$.
\end{proof}

\begin{corr}
	When $n\geq2$, the equivariant harmonic maps from $\Cn$ to $\CPn$ all have infinite energy.
\end{corr}

\section{Self-similar solutions}
\label{sec:selfsimilar}

In this section we study the
self-similar solutions, 
which are solutions of the form $u(r,t) = \psi(r/\sqrt{t})$ 
for a profile $\psi(r) = u(r,1)$.

To determine a convenient equation for the profile, we take the \gll{} flow PDE \eqref{eqn:2:sphere} and multiply both sides by $(\alpha u \times + \beta P)$.
Using the relationship,
\[
	(\alpha u \times + \beta P ) (\alpha P + \beta u\times) = (\alpha^2+\beta^2) u\times = u \times,
\]
(compare to $(\alpha i +\beta)(\alpha+\beta i )=i$)
we may equivalently write the PDE as,
\begin{equation}
	\label{eqn:3:sphere2}
	\alpha u\times u_t + \beta u_t = u \times \left( \der{^2 u }{r^2} + \frac{2n-1}{r} \der{u}{r}  + \frac{ 2n-2+u_3 }{r^2} e_3 \right).
\end{equation}
We now substitute in $u(r,t) = \psi(r/\sqrt{t})$ to determine the ODE for the profile.

\begin{defn}
	The self-similar problem for the \gll{} flow is given by the ODE,
	\begin{equation}
		\label{eqn:3:selfsim}
		-\frac{r}{2} \left( \alpha \psi \times \psi_r + \beta \psi_r \right)= \psi\times\left( \sder{\psi}{r} + \frac{2n-1}{r} \der{\psi}{r} + \frac{(2n-2+\psi_3)}{ r^2 } e_3 \right) 
	\end{equation}
	subject to the initial conditions $\psi(0)=0$ and $\psi'(0)=v = (v_1,v_2,0) \in T_{e_3} \Sp^2 $.
\end{defn}

In the following sequence of Lemmas we will prove Theorem \ref{thm:1}, as stated on page \pageref{thm:1} in the introduction.

\begin{lemma}
	For every $v=(v_1,v_2,0) \in T_{e_3} \Sp^2$ there is a unique global solution to \eqref{eqn:3:selfsim}.
	For $r>0$ this global solution is smooth and, if $v\neq 0$, satisfies $\psi(r)\neq e_3$.
\end{lemma}

\begin{proof}
Local existence and uniqueness in a neighborhood of the singular point $r=0$ 
follows from the Theorem \ref{thm:app:wellposedness} in the Appendix.
For $r>0$, the ODE \eqref{eqn:3:selfsim} is smooth and local existence, uniqueness and smoothness comes from the standard ODE theory.
In order to prove global existence we establish an \emph{a priori} bound on the derivative of $\psi$.

Define the function $A(r) = r^2 |\psi_r|^2$
We have,
\begin{equation}
	\label{eqn:3:Aprime}
	A'(r) = 2 r |\psi_r|^2 + 2 r^2 \psi_{rr} \cdot \psi_r.
\end{equation}
In order to calculate $\psi_{rr} \cdot \psi_r$, we take the inner product of the ODE \eqref{eqn:3:selfsim} with $\psi \times \psi_r$.
Using the fact that if $v$ or $w$ is orthogonal to $u$, then $(u\times v)\cdot(u\times w) = v\cdot w$,
and also the relation $v \cdot(u\times v)=0$, we determine that,
\[
	-\frac{\alpha r}{2} |\psi_r|^2 = \psi_{rr} \cdot \psi_r + \frac{2n-1}{r} | \psi_r |^2 + \frac{ 2n-2+\psi_3 }{ r^2 } e_3 \cdot \psi_r,
\]
and hence by solving for $\psi_{rr}\cdot \psi_r$ and substituting this into \eqref{eqn:3:Aprime} we find,
\begin{align}
	A'(r) &= 2 r |\psi_r|^2 -\left( \frac{2n-1}{r} + \frac{\alpha r}{2} \right) 2 r^2 |\psi_r|^2  - (2n-2+\psi_3)(\psi_3)_r	\nonumber \\
		&= -\left( \frac{2n-3}{r} + \frac{\alpha r}{2} \right) 2 A(r) - \frac{d}{dr} \left[ (2n-2) \psi_3 + \frac{\psi_3^2}{2} \right]. \label{eqn:3:aprioriode}
\end{align}
Integrating this equation gives,
\begin{equation}
	\label{eqn:3:apriori}
	A(r) + \int_0^r \left( \frac{2n-3}{s} + \frac{\alpha s}{2} \right) 2 A(s)ds = \left[ (2n-2) (1-\psi_3) + \frac{1-(\psi_3)^2}{2} \right].
\end{equation}
To bound $A(r)$, we observe that the integral on the left hand side is non-negative because $A(s) \geq 0$,
	and so the left hand side is bounded below by $A(r)$.
On the other hand, we have $\psi_3 \in [-1,1]$ and hence the right hand side is bounded above by $4n$.
This then gives $A(r) \leq 4n$, and $|\psi_r| \leq 2n/r$.
This proves global existence.
(The constants $4n$ and $2n$ are, of course, not optimal; 
	they are noted merely to 
        show that the constants may be chosen independently of $\psi$.)

To prove that $\psi(r) \neq e_3$ for $r>0$ we observe that the integral on the left hand side in \eqref{eqn:3:apriori} is increasing in $r$.
In the non-trivial case $v\neq 0$, it is strictly increasing a neighborhood of $r=0$ because 
$A'(r) = r^2 |\psi_r|^2 \geq \epsilon r^2$ in a neighborhood of $r=0$.
Hence in this case the integral is strictly positive for $r>0$.
Because $A(r) \geq 0$ we see that the left hand side of \eqref{eqn:3:apriori} is strictly positive and so,
\[
	\left[ (2n-2) (1-\psi_3) + \frac{1-(\psi_3)^2}{2} \right] > 0, 
\]
for $r>0$.
This gives  $\psi_3(r) \neq 1$, which means $\psi(r) \neq e_3$.
\end{proof}

\begin{lemma}
	If $\alpha >0$ we have $|\psi_r| \lesssim 1/r^3$.
\end{lemma}

\begin{proof}
	Recall the bound $A(r) \leq 4n$.
	Using equation \eqref{eqn:3:aprioriode} we have
	\begin{align*}
		A'(r)	&\leq -\frac{\alpha r }{2} A(r) - (2n-2+\psi_3)(\psi_3)_r	
                        \leq -\frac{\alpha r }{2} A(r) + \left| \frac{2n-2+\psi_3}{ r^{3/2} \sqrt{\alpha/2} } \right| \cdot \left| r^{3/2} \sqrt{\alpha/2} (\psi_3)_r \right|	\\
			&\leq -\frac{\alpha r }{2} A(r) + \frac{1}{2} \left( \frac{ 8 n^3 }{ r^3 \alpha } +  \frac{ r^3 \alpha }{2} | \psi_r|^2  \right)		
			= -\frac{\alpha r }{4} A(r) + \frac{ 4 n^3 }{ r^3 \alpha }.
	\end{align*}
	Integrating this equation then gives $A(r) \lesssim A(1)e^{-\alpha r^2/8} + 1/r^4 \lesssim 1/r^4$
		and $|\psi_r| \lesssim 1/r^3$.
	(The details of how this integration may be performed are given in Proposition \ref{prop:app:integration} in the appendix.)
\end{proof}

\begin{lemma}
	There exists a point $\psi_\infty \in \Sp^2$, $\psi_\infty \neq e_3$, such that $\lim_{r \rightarrow \infty} \psi(r) = \psi_\infty$.
	We have the convergence rate inequality $|\psi_\infty-\psi(r)|\leq 40n^2/r^2$.
	The profile $\psi$ has infinite energy.
\end{lemma}

\begin{proof}
	For $\alpha>0$, the bound $|\psi_r|\lesssim 1/r^3$ implies convergence of $\psi$ in the limit $r \rightarrow \infty$.
	In the case $\alpha=0$, when there is no heat flow contribution, the decay on the derivative is less strong,
		and so a different argument is needed.
	However in the proof we consider the general case as it is useful to know that the constant in the rate of convergence equation
		may be chosen independently of $\psi$.

	We first multiply the ODE
        \eqref{eqn:3:sphere2} by $(-\alpha \psi \times + \beta P )$.
	We have the relations  $(-\alpha \psi \times + \beta P )( \alpha \psi \times + \beta P)=(\alpha^2+\beta^2) P = P$
		and $(\psi\times)(\psi \times) =-P$ (compare to $(-\alpha i + \beta)(\alpha i + \beta) = 1$ and $(i)(i)=-1$).
	We can thus write the equation as,
\begin{align*}
	-\frac{ r }{2} \psi_r 
		&= (\alpha P + \beta \psi \times) \left(  \psi_{rr} + \frac{2n-1}{r} \psi_{r}  + \frac{ 2n-2+\psi_3}{ r^2 } e_3 \right), \\
		&= (\alpha  + \beta \psi \times) \left( \frac{1}{r^{2n-1}} \der{}{r} ( r^{2n-1} \psi_r ) + |\psi_r|^2 \psi  + \frac{ 2n-2+\psi_3}{ r^2 } P e_3 \right),
\end{align*}
where in the second equality we have moved the projection $P$ inside and  expanded $P \psi_{rr} = \psi_{rr} + |\psi_r|^2 \psi$.
We divide through by $r$ and integrate over $[r_1,r_2]$ to determine that,
\begin{align*}
 	-\frac{1}{2} ( \psi(r_2) - \psi(r_1 ) )	
		&=
			\int_{r_1}^{r_2} (\alpha + \beta \psi \times) \left( \frac{1}{r^{2n}} \der{}{r} ( r^{2n-1} \psi_r ) + \frac{ |\psi_r|^2 \psi }{r} + \frac{ 2n-2+\psi_3}{ r^3 } P e_3 \right) dr.
\end{align*}
Now integrating by parts in the first term yields,
\begin{align*}
 	-\frac{1}{2} ( \psi(r_2) - \psi(r_1 ) )	
		&= 
		\frac{ [\alpha + \beta \psi(r_2) \times ]\psi_r(r_2) }{ r_2 }
		-\frac{ [\alpha + \beta \psi(r_1) \times ] \psi_r(r_1) }{ r_1 }	\\
		&\hspace{1cm} -\int_{r_1}^{r_2}  (\alpha + \beta \psi \times) \left( \frac{-2n}{r^{2n+1} } r^{2n-1} \psi_r \right) dr	\\
		&\hspace{1cm} +  \int_{r_1}^{r_2} \left( \alpha \frac{|\psi_r|^2 \psi}{r} + \frac{2n-2+u_3}{r^3} (\alpha P + \beta  \psi \times ) e_3 \right) dr.
\end{align*}
Now using the bounds $|\psi(r)|=1$ and $|\psi_r(r)|\leq 2n/r$ yields,
\[
	\frac{1}{2} | \psi(r_2) - \psi(r_1) | \leq \frac{2n}{r_2^2} + \frac{2n}{r_1^2} + \int_{r_1}^{r_2} \frac{4n^2}{r^3} dr
		 + \int_{r_1}^{r_2}\left( \alpha \frac{4n^2}{r^3} + \frac{2n}{r^3} \right) dr \leq \frac{20 n^2}{r_1^2},
\]
which implies the solution converges with the rate given in the statement of the Lemma.

To see that the limit $\psi_\infty$ cannot be $e_3$ we consider equation \eqref{eqn:3:apriori} again.
As discussed previously, the integral in \eqref{eqn:3:apriori} is strictly positive and non-increasing for $r>0$.
If $\delta >0$ denotes the value of the integral at $r=1$ we then have, for all $r>1$,
\[
	\delta \leq  \int_0^r \left( \frac{2n-3}{s} + \frac{\alpha s}{2} \right) 2 A(s)ds \leq \left[ (2n-2) (1-\psi_3(r)) + \frac{1-(\psi_3(r))^2}{2} \right].
\]
We therefore have
\[
	\delta \leq  \left[ (2n-2) (1-\psi_3(\infty)) + \frac{1-(\psi_3(\infty))^2}{2} \right],
\]
which gives $\psi_\infty \neq e_3$.

Because the  limit is not $e_3$, the profile has infinite energy by Lemma \ref{lemma:2:infenergy}.
\end{proof}

\begin{lemma}
	When $\alpha=0$ we have $\lim_{r \rightarrow\infty} r|\psi_r| =0$.
\end{lemma}

\begin{proof}
	It is sufficient to show that $\lim_{r \rightarrow \infty} A(r) =0 $.
	In the $\alpha=0$ case equation \eqref{eqn:3:apriori} reads.
	\[
		A(r) + \int_0^r \left( \frac{2n-3}{s} \right) 2 A(s)ds = \left[ (2n-2) (1-\psi_3) + \frac{1-(\psi_3)^2}{2} \right].
	\]
	We know from the previous lemma that $\psi_3$ converges as $r \rightarrow \infty$.
	The integral also converges simply because it is non-decreasing; moreover, because it is bounded above (by $4n$) it converges to a real number.
	We then have that $A(r)$ converges as $r\rightarrow \infty$.
	By examining the integral, which is finite in the limit, we see that we must have $\lim_{r \rightarrow \infty} A(r)=0$.
\end{proof}

\begin{lemma}
	The limit $\psi_\infty$ is a continuous function of the initial data $v$.
	In particular, as $v \rightarrow 0$ we have $\psi_\infty \rightarrow 0$.
\end{lemma}

\begin{proof}
	For convenience we will denote the self-similar profile corresponding to initial data $v$ by $\psi_v(r)$,
		and we will let $\psi_v(\infty)$  denote its limit as $r \rightarrow \infty$.

	The ODE local existence results give that for any $r_0>0$ the map $v \mapsto \psi_v(r_0)$ is continuous.

	We have previously established the bound, for $r_1<r_2$,
	$
		|\psi_v(r_2) - \psi_v(r_1)| \leq { 60 n^2  }/{ r_1^2 }
	$
	This shows that the map $v \mapsto \psi_v(r)$ converges to the map $v \mapsto \psi_v(\infty)$ uniformly,
		and hence that the map $v \mapsto \psi_v(\infty)$ is continuous.
	
	Finally, we note that $\psi_0(r) \equiv 0$, $\psi_0(\infty)=0$,
		and so $\lim_{v \rightarrow 0 } \psi_v(\infty) =0$, by continuity.
\end{proof}

With this Lemma, the proof of Theorem \ref{thm:1} is complete.

\section{Global critical wellposedness in dimension two}
\label{sec:wellposedness}

In this section we prove a global critical small data wellposedness theorem
for the \Sch{} maps equation for
	equivariant maps from $\C^n$ to $\CPn$ when $n=2$.
The equation may be written in the sphere coordinates as,
	\begin{equation}
		\label{eqn:4:sphere}
		u_t(r,t) = u \times \left( \der{^2 u}{r^2} + \frac{2n-1}{r} \der{u}{r} + \frac{ 2n-2+u_3 }{r^2} e_3  \right), 	
	\end{equation}
	or equivalently as,
	\begin{equation}
		\label{eqn:4:sphere2}
		- u \times u_t(r,t) = \der{^2 u}{r^2}  + |u_r|^2 u + \frac{2n-1}{r} \der{u}{r} + \frac{ 2n-2+u_3 }{r^2} P_u e_3 ,
	\end{equation}
        where $P_u e_3$ is the projection of the vector $e_3=(0,0,1)$
        onto the tangent space at $u$.

Our proof relies on techniques that have been developed for the 
\Sch{} maps equation for the sphere.
Because of the structural 
similarity between that equation and \eqref{eqn:4:sphere2},
such techniques can be adapted here.
We first use a form of the Hasimoto transform to determine
an equation on a derivative term of $u$ that has a
simpler nonlinearity.
We then formulate the fixed point argument, and determine
necessary estimates on the nonlinearity for the fixed point
argument to be carried through.
We conclude by proving these estimates in the case $n=2$,
    thereby establishing Theorem 2.

We present our work in terms of the 
\Sch{} maps equation $(\alpha=0)$, however our proof is 
valid for the general
\gll{} case when $\beta>0$ because all the same estimates
(in particular the Strichartz estimates)
still apply. 

\subsection{Derivation of the PDE through the Hasimoto transform}

The Hasimoto transform is an extensively used tool for proving wellposedness 
of \Sch{} maps equations when the target is the sphere
	or a general complex surface.
In geometric terms, it arises as follows.
For fixed $t$, a smooth solution of \eqref{eqn:4:sphere} will satisfy $u(0,t)= e_3$.
The function $r \mapsto u(r,t)$ thus defines a curve in $\Sp^2$ starting
at $e_3$ at $r=0$.
If one fixes a unit tangent vector $e(0) \in T_{e_3} \Sp^2$, one can consider
the parallel transport $e(r)$ of this 
	vector along the curve $r \mapsto u(r,t)$; the function $e(r)$
        satisfies $D_r e(r) = \nabla_{u_r}e(r) = 0$.
Now because the tangent space at any point is two dimensional, 
	the vectors $e(r)$ and $Je(r)$ give a basis for the tangent space
        $T_{u(r)} \Sp^2$.
Any derivative of $u$, or other element of the tangent space, can be expressed in terms of this basis.
In our case, we define a complex valued function $q$ by the formula,
\begin{equation}
	\label{eqn:4:hasi}
	qe = (\Re q  + \Im q \, J)e =  u_r.
\end{equation}
We then determine an equation on $q$.
The right hand side is chosen so that $q$ will satisfy a \Sch{} equation
with a non-linearity that is easier to handle than that of 
\eqref{eqn:4:sphere2}.

\begin{lemma}
    The function $q$ satisfies the PDE,
    \begin{equation}
        i q_t = q_{rr} + \frac{2n-1}{r} q_r - \frac{2n-1}{r^2} q + N(q),
        \label{eqn:4:qpde1}
    \end{equation}
    where the nonlinear term $N(q)$ is given by,
    \begin{equation}
            N(q) =  \frac{d}{dr} \left[
        - \frac{2n-2 + u_3 }{ r^2 } \int_0^r u_3(s) q(s) ds
            \right] + \alpha q,
        \label{eqn:4:qpde2}
    \end{equation}
    for a real-valued function $\alpha$ satisfying,
    \begin{equation}
        \alpha_r =  \Re \left( 
            \overline{q} q_r + \frac{|q|^2}{r} 
            - \overline{q} \frac{2n-2 + u_3 }{ r^2 } \int_0^r u_3(s) q(s) ds
        \right).
        \label{eqn:4:qpde3}
    \end{equation}
\end{lemma}

\begin{proof}
    First, we recall that in the embedding $\Sp^2 \subset \R^3$
    the covariant derivative of a vector field $v(r) \in T_{u(r)}\Sp^2$
    is given by $D_r v = v_r + \langle u_r, v \rangle u$, where the inner product
    here is the usual inner product on $\R^3$.

    Now let $p$ and $q$ satisfy $pe = u_t$ and $qe = u_r$.
    We will determine three equations relating $p$, $q$ and $u$.

    \begin{enumerate}
    \item
    Because $e$ satisfies $D_re=0$ we have,
    \begin{equation}
        q_re = D_r (qe) = D_r( u_r) = u_{rr} + |u_r|^2u,
        \label{eqn:4:urr}
    \end{equation}
    which is the first two terms two term in the right hand side of \eqref{eqn:4:sphere2}.
    The next term in \eqref{eqn:4:sphere2} is $((2n-1)/r)qe$.
    For the projection term we calculate, using $D_re=0$,
    \begin{align*}
        \frac{d}{dr} \langle P_u e_3, e \rangle
        &= \frac{d}{dr} \langle e_3 - \langle u, e_3 \rangle u , e \rangle
        = \langle D_r( e_3 - \langle u, e_3 \rangle u ), e \rangle  \\
        &= \left\langle 
            \frac{d}{dr}( e_3 - \langle u, e_3 \rangle u )
            + \langle u_r, e_3 - \langle u, e_3 \rangle u \rangle u,
            e \right\rangle \\
        &= \left\langle 
            - \langle u_r, e_3 \rangle u 
            - \langle u, e_3 \rangle u_r
            + \langle u_r, e_3  \rangle u,
            e \right\rangle  = - u_3 \langle u_r ,e \rangle 
            = - u_3(r) \Re p(r).
    \end{align*}
    Using the fact that $u(0,t) = e_3$, so that $P_u e_3=0$ at $r=0$, we have,
    \begin{equation}
        \langle P_u e_3, e \rangle = - \int_0^r u_3(s) \Re q(s) ds.
        \label{eqn:4:Pucoords}
    \end{equation}
    An identical calculation for $\langle P_u e_3, Je \rangle$ gives,
    in total,
    \[
        P_u e_3 = 
        - \left( \int_0^r u_3(s) q(s) ds \right) e(r).
    \]

    Plugging \eqref{eqn:4:urr} and \eqref{eqn:4:Pucoords} into 
    \eqref{eqn:4:sphere2} then gives,
    \begin{equation}
        ip = q_r + \frac{2n-1}{r} q 
            - \frac{2n-2+u_3}{r^2} \int_0^r u_3(s) q(s) ds.
            \label{eqn:4:ip}
    \end{equation}

    \item
    From the identity $D_r u_t = D_t u_r$ we find,
    \begin{equation}
        p_r e = D_r( pe ) = D_r u_t = D_t u_r = D_t(qe) = q_t e + q D_te.
        \label{eqn:4:comm}
    \end{equation}
    Because $e$ is a parallel transport vector field, $|e|^2=1$ and so
    $ 0 = (d/dt) |e|^2 = \langle D_t e ,e \rangle$.
    The vector $D_t e$ is thus orthogonal to $e$.
    Because the tangent space is spanned by $e$ and $Je$,
    we must have $D_te = \alpha Je$ for some real-valued function $\alpha$.
    Substituting this into \eqref{eqn:4:comm}, we get
    $p_r e = q_t e + q \alpha Je$, or,
    \begin{equation}
        p_r = q_t + i \alpha q.
        \label{eqn:4:pr2}
    \end{equation}

    \item
    To determine an equation on $\alpha$ we use the curvature relation
    $D_t D_r e = D_r D_t e + R(u_t, u_t) e$
    where $R$ is the Riemann curvature tensor.
    On the sphere $R(v,w)z = \langle J v, w \rangle Jz$.
    Therefore, using also $D_re=0$, we find,
    \[
        0 = D_r(\alpha J e) + \langle J u_t, u_r \rangle Je
            =   \alpha_r Je +  \langle p J e, qe \rangle Je,
    \]
    which gives $\alpha_r = - \Im (p \overline{q})$.
    Substituting the formula for $p$ in \eqref{eqn:4:ip} gives 
    equation \eqref{eqn:4:qpde3}.
    \end{enumerate}

    To determine an equation only on $q$ we differentiate 
    \eqref{eqn:4:ip} with respect
    to $r$, to find,
    \[
        ip_r = q_{rr} + \frac{2n-1}{r} q_r - \frac{2n-1}{r^2} q
            +  \frac{d}{dr} \left[
        - \frac{2n-2 + u_3 }{ r^2 } \int_0^r u_3(s) q(s) ds,
            \right].
    \]
    Substituting the expression for $p_r$ in
    \eqref{eqn:4:pr2} gives equations 
    \eqref{eqn:4:qpde1} and \eqref{eqn:4:qpde2}.
\end{proof} 

\subsection{Formulating the fixed point argument}

We recall Theorem \ref{thm:2} from the introduction.

\begin{thm*}[Theorem \ref{thm:2}, page \pageref{thm:2}]
	Fix $p \in [1,2]$ and define
	\[
		 \frac{1}{r} = \frac{1}{2} - \frac{1}{6p}
	\]
	and the spaces $X$ and $X_0$ given by the norms
	\[
		\| q \|_X = \| \nabla q \|_{L^{3p}_t L^{r}_x}
		\;\;\text{ and }\;\;
		\| q \|_{X_0} = \| e^{it\Delta} (aq) \|_{X},
	\]
	where $a(x) = x_1/r$.
        There exists $\epsilon>0$ such that if 
            $\| q_0 \|_{X_0} \leq \epsilon$
        there is unique global solution of \eqref{eqn:4:qpde1} 
            for $n=2$ in the space $X$.
\end{thm*}

We begin by determining a convenient Duhamel representation for the problem.
Our Duhamel representation will be valid for all $n$,
though we carry out the wellposedness argument for $n=2$ only.
In the following we will rely heavily 
on the Hardy inequalities given in Theorems \ref{thm:app:hardy1}
and \ref{thm:app:hardy2} in the appendix.

First, we absorb the linear term $-(2n-1)q/r^2$ into the Laplacian.
To do this, we fix a function $a:\Sp^{2n-1} \rightarrow \C$ that
satisfies $\Delta_{\Sp^{2n-1}}a = -(2n-1) a$.
We may concretely choose $a(x) = x_1 $.
To see this, extend $a$ to a function on all of $\R^{2n}$ by $a(x/|x|)$.
On the one hand, we have,
\[
	\Delta_{\R^{2n}} \left( r a (x/|x|) \right) 
            = \Delta_{\R^{2n}} \left( x_1 \right) = 0.
\]	
Then, using the polar representation,
$\Delta_{\R^{2n}} = \partial_{rr} + ((2n-1)/r) \partial_r + (1/r^2) \Delta_{\Sp^{2n-1}}$,
we see that,
\begin{align*}
	0 =
        \left[ \partial_{rr} + \frac{2n-1}{r} \partial_r 
        + \frac{1}{r^2} \Delta_{ \Sp^{2n-1}} \right] (r a(x/|x|))
        &=	\left[ 0 + \frac{2n-1}{r}  \right] a(x/|x|) 
        + \frac{r}{r^2} \Delta_{\Sp^{2n-1}}  a(x/|x|) ,
\end{align*}
and so,
\[
	\Delta_{\Sp^{2n-1}} a(x/|x|) = -(2n-1) a(x/|x|).
\]
Now defining $w(x,t) = q(r,t) a(x/|x|)$, we see that,
\begin{equation}
	\label{eqn:4:lapw}
	\Delta_{\R^{2n}} w = \der{^2 q}{r^2} a + \frac{2n-1}{r} \der{q}{r} a -\frac{2n-1}{r^2} qa.
\end{equation}
This is exactly the Laplacian term in the PDE \eqref{eqn:4:qpde1}
multiplied by $a$.

In terms of estimates, we have the pointwise estimate $|\nabla a| \leq 1/r$,
which is determined from a calculation.
For Lebesque estimates we have,
\begin{align*}
	\| w \|_{L^p}^p 
		= \| qa \|_{L^p}^p &= \int_0^\infty \left( |q(r)|^p \int_{r \Sp^{2n-1} } | a(x/|x|)|^p dx \right) dr		\\
		&= \int_0^\infty \left( |q(r)|^p r^{2n-1} \int_{\Sp^{2n-1} } |a(x/|x|)|^p dx \right) dr		\\
		&= C \| q \|_{L^p}^p,
\end{align*}
where $C = \| a \|_{L^p(\Sp^{2n-1})}/|\Sp^{2n-1}| < \infty$.
We also have,
\[
        \| \nabla q \|_{L^p} \sim \| \nabla_r q \|_{L^p}
        \sim 
            \| \nabla_r(aq) \|_{L^p}
        \lesssim 
            \| \nabla_r(aq) \|_{L^p}
            + \| \nabla_\theta(aq) \|_{L^p} \sim \| \nabla w \|_{L^p},
\]
while,
\[
    \| \nabla w \|_{L^p} \sim
        \| \nabla(a) q \|_{L^p} +
        \| a \nabla_r q \|_{L^p} 
        \lesssim
        \| \frac{1}{r} q \|_{L^p} +
        \| \nabla_r q \|_{L^p} 
        \lesssim \| \nabla q \|_{L^p},
\]
so in conclusion $\| \nabla q \|_{L^p} \sim \| \nabla w \|_{L^p}$.

In the next Lemma we rather carefully verify that we can recover solutions to the PDE for $q$ from solutions to the PDE for $w$.

\begin{lemma}
	The PDE on $w$ is given by,
	\begin{equation}
		\label{eqn:4:wpde}
		w_t  = \Delta w + N(q) a
	\end{equation}
	or in Duhamel form by,
	\begin{equation}
		w(x,t) = e^{it\Delta} w(x,0) + i \int_0^r e^{i(t-s)\Delta} N(q(r,s)) 
                    a(x/|x|,s) ds
                    \label{eqn:4:wduh}
	\end{equation}

	If the solution $w(x,t)$ corresponding to initial data of the form $w(x,0) = q(r)a(x/|x|)$ is unique,
		then the solution is of the form $w(x,t) = q(r,t) a(x/|x|)$,
		where $q$ satisfies \eqref{eqn:4:qpde1}.
\end{lemma}

\begin{proof}
	To determine the equation \eqref{eqn:4:wpde}
        for $w$ we simply multiply the PDE for $q$ \eqref{eqn:4:qpde1}
        by $a$, and use the expression \eqref{eqn:4:lapw} for $\Delta w$.
	The Duhamel representation is standard.

	We now show how solutions of \eqref{eqn:4:qpde1} may be recovered 
        from solutions of the equation for $w$.
        Let $w$ be a solution of \eqref{eqn:4:wduh}
	and define $\tilde{w} = -(1/(2n-1)) \Delta_{\Sp^{2n-1}}w$.
	Assuming uniqueness we will show that $\tilde{w}=w$.
	We take the spherical Laplacian $-(1/(2n-1))\Delta_{\Sp^{2n-1}}$ of \eqref{eqn:4:wpde}, noting that it commutes both with $\Delta=\Delta_{\R^{2n}}$ 
        and $N(q)$, as $N(q)$ is radial.
	We then find that $\tilde{w}$ satisfies the same PDE \eqref{eqn:4:wpde} as $w$.
	Moreover, we have,
	\[
		\tilde{w}(x,0) = - \frac{1}{2n-1} \Delta_{\Sp^{2n-1}} w(x,0) =
			- \frac{1}{2n-1} \Delta_{\Sp^{2n-1}} \left[ q(r) a(x/|x|)\right]  
                        = w (x,0),
	\]
	and so by uniqueness, $\tilde{w}(x,t) = - (1/(2n-1)) \Delta_{\Sp^{2n-1}}w(x,t) = w(x,t)$.
	This means that $w$ is a radial function times an eigenfunction of the Laplacian of the sphere of $\Sp^{2n-1}$ with eigenvalue $-(2n-1)$.	

	Let $T_k:\R^{2n} \rightarrow \R^{2n}$ be the linear map that multiplies the $k$th component 
		of $x \in \R^{2n}$ by $-1$ and leaves the other components fixed.
	From the representation of $a$ we see that for $k=1$ we have $w_0(T_k x) = - w_0(x)$
		while for $k\geq 2$ we have $w_0(T_k x) = w_0(x)$.
	By uniqueness, $x\mapsto w(x,t)$ inherits these properties also.
	But now the only eigenfunction of the Laplacian on the sphere with eigenvalue
            $-(2n-1)$ with these symmetries is precisely $a$.
	Therefore $w(x,t) = q(r,t) a(x/|x|)$.
	Substituting this expression
        into the PDE \eqref{eqn:4:wpde} for $w$ yields the
        PDE \eqref{eqn:4:qpde1} for $q$.
\end{proof}

By virtue of this Lemma, we can perform the fixed point argument on $w$.
The next Lemma describes sufficient estimates for this fixed point argument to hold,
	and in the proof the fixed point argument is described.

\subsection{Wellposedness when $n=2$}

For the reminder of this section we fix $n=2$.

Before stating the lemma we fix some index notation.
In the course of the proof we will need to handle Lebesgue space
norms of quantities like $q$, $q_r$, $q^2$, $q q_r$, etc., and other quantities which scale like these.
We are led to define the index,
\begin{equation}
	\label{eqn:4:s}
	\frac{1}{ s(i,j) } = \frac{ i + j }{ 4 } - \frac{ i }{ 6p }.
\end{equation}
We will put items that scale like the product of $i$ copies of $q$ with a total of $j$ derivatives  in the space $\LebS{s(i,j)}$.
For example, we will put $q$ in $\LebS{s(1,0)}$, we will put $q^2$ in $\LebS{s(2,0)}$ and $qq_r$ in $\LebS{s(2,1)}$.
In this way, critical scaling is maintained throughout as, for example, $\|qq_r \|_{\LebS{s(2,1)}}$ is invariant under scaling.

The Strichartz inequality we will use is,
\[
	\left\| \int_0^t e^{i(t-s)\Delta} G ds \right\|_{\LebST{3p}{r}} \lesssim \left\| G \right\|_{\LebST{p}{s(3,1)}};
\]
this is classical: see, for example, \cite{Tao}.
The \Hol{} inequality is,
\[
	\| f g \|_{\LebS{s(i+k,j+m)}} \leq \| f \|_{\LebS{s(i,j)}} \cdot \| g \|_{\LebS{s(k,m)}};
\]
and the Sobolev is, for $k<l$,
\[
	\| \nabla^k f \|_{\LebS{s(i,j)}} \lesssim \| \nabla^l f \|_{\LebS{ s(i,j+l-k) }}.
\]
One verifies that these inequalities hold by checking the relevant
exponent conditions.

Finally, note that $s(1,1)=r$.

\begin{lemma}
	For Theorem 2 to be true, it is sufficient that the following bounds hold;
	\begin{align}
		\label{eqn:4:Nbounds}
		\| \nabla N(q) \|_\LebS{s(3,1)} &\lesssim \|\nabla q \|_\LebS{s(1,1)}^3 \\ 
		\| \nabla ( N(q_1) - N(q_2)) \|_\LebS{s(3,1)} &
                \lesssim \|\nabla (q_1 - q_2) \|_\LebS{s(1,1)} \left(\|\nabla q_1 \|_{\LebS{s(1,1)}}^2 + \| \nabla q_2 \|_\LebS{s(1,1)}^2	\right)
		\label{eqn:4:Nbounds2}
	\end{align}
\end{lemma}
\begin{proof}
Well-posedness follows by a fixed point argument for  the operator,
\[
	T w = e^{it\Delta} w(x,0) + i \int_0^t e^{i(t-s)\Delta} N(q(r,s)) a(x,s) ds.
\]	
We will show that $T$ is a contraction mapping on a small ball around 0.

We first show that $T$ maps a ball to itself.
We have the bound,
\begin{align*}
	\| Tw \|_X
		&\leq \| w_0 \|_{X_0}  + \left\| \int_0^t e^{i(t-s) \Delta} \nabla( N(q) a ) ds \right\|_{\LebST{3p}{r}}	
		\lesssim \| w_0 \|_{X_0} + \left\| \nabla( N(q) a )  \right\|_{\LebST{p}{s(3,1)}}.
\end{align*}
Considering the space norm of the integral, we have, by \Hol{} and Sobolev, and then conditions \eqref{eqn:4:Nbounds},
\begin{align*}
	\| \nabla ( N(q) a) \|_{\LebS{s(3,1)}}
		&\lesssim \| \nabla(N) a \|_{\LebS{s(3,1)}} 
		    + \| N \nabla a \|_{\LebS{s(3,1)}}  \\
		&\leq \| \nabla(N) \|_{\LebS{s(3,1)}} 
		    + \left\|  \frac{1}{r} N  \right\|_{\LebS{s(3,1)}} 
		\leq \| \nabla(N) \|_{\LebS{s(3,1)}} 
		\lesssim \| \nabla q \|_{\LebS{s(1,1)}}^2  
		\lesssim \| \nabla w \|_{\LebS{s(1,1)}}^3,
\end{align*}
and hence, as $r=s(1,1)$,
\[
	\| Tw \|_X \leq \| w_0 \|_{X_0} + C \| w \|_X^3.
\]
Now choose $\epsilon_0$ so that $C \epsilon_0^2 \leq 1/2$,
and let $\epsilon \leq \epsilon_0$.
Then, if $\| w_0 \|_{X_0} \leq \epsilon/2$ and $\|w\|_X \leq \epsilon$ we have,
\[
	\| Tw \|_X \leq \frac{\epsilon}{2} + C \epsilon^3 \leq \frac{\epsilon}{2} + (C\epsilon_0^2) \epsilon \leq \epsilon,
\]
and so $T$ maps every $\epsilon$ ball into itself, for $\epsilon$ sufficiently small, 
assuming the initial data satisfies the bound $\|w_0\|_{X_0} \leq \epsilon/2$.

We next show that $T$ is a contraction in a sufficiently small ball around 0.
Let $w_1$ and $w_2$ be two solutions, with radial 
parts $q_1$ and $q_2$ respectively. We have,
\[
	T w_1 - Tw_2 = \int_0^t e^{i(t-s)\Delta} ( \nabla(N(q_1)a) - \nabla(N(q_2)a)) ds ,
\]
which gives, using \eqref{eqn:4:Nbounds2},
\begin{align*}
	&\hspace{-1cm}\left\| \nabla(N(q_1)a) - \nabla(N(q_2)a) \right\|_{\LebS{s(3,1)}}	\\
        &\lesssim 
        \left\| \nabla(N(q_1)-N(q_2)) a 
            \right\|_{\LebS{s(3,1)}}	+ 	
        \left\| ( N(q_1)-N(q_2) )\nabla a 
            \right\|_{\LebS{s(3,1)}}		
        \\
        &\lesssim 
        \left\| \nabla(N(q_1)-N(q_2)) a 
            \right\|_{\LebS{s(3,1)}}	+ 	
        \left\|  \frac{1}{r} ( N(q_1)-N(q_2)) 
            \right\|_{\LebS{s(3,1)}}		
        \\
	& \lesssim (\| \nabla q_1 \|_{\LebS{s(1,1)}}^2 + \| \nabla q_2 \|_\LebS{s(1,1)}^2 ) \| \nabla (q_1 - q_2) \|_\LebS{s(1,1)}, \\
	& \lesssim (\| \nabla w_1 \|_{\LebS{s(1,1)}}^2 + \| \nabla w_2 \|_\LebS{s(1,1)}^2 ) \| \nabla (w_1 - w_2) \|_\LebS{s(1,1)},
\end{align*}
and so,
\[
	\| Tw_1 - Tw_2 \|_X \lesssim
		(\| w_1 \|_X^2 + \| w_2 \|_X^2) \| w_1 - w_2 \|_X
\]
and hence by choosing the ball small enough, $T$ is a contraction.
\end{proof}

\begin{lemma}
When $n=2$ the bounds \eqref{eqn:4:Nbounds} and \eqref{eqn:4:Nbounds2} hold.
\end{lemma}

\begin{proof}
Write,
\[
    N = \frac{d}{dr} \left(
        - \frac{2n-2 + u_3 }{ r^2 } \int_0^r u_3(s) q(s) ds
        \right) + \alpha q =: N_1 + N_2,
\]
and recall,
\begin{equation}
    \alpha_r = \Re \left( q_r \overline{q} + \frac{|q|^2}{r} -
        \overline{q}  \frac{2n-2 + u_3 }{ r^2 } \int_0^r u_3(s) q(s) ds
        \right).
        \label{eqn:4:alphar}
\end{equation}
We will prove the bounds for $N_1$ first.

We have,
\begin{align*}
    \| \nabla N_1 \|_{s(3,1)}
    &\lesssim
    \left\| \nabla^2  \left(
        - \frac{2n-2 + u_3 }{ r^2 } \int_0^r u_3(s) q(s) ds
        \right)
        \right\|_{\LebS{s(3,1)}}   \\
    &\lesssim
    \left\|
        (\nabla^2 u_3) \frac{1}{r^2} \int_0^r u_3(s)q(s)ds 
    \right\|_{\LebS{s(3,1)}} \\
    &\hspace{1cm} + \left\|
        - \frac{2n-2 + u_3 }{ r^2 } \frac{1}{r^2} \int_0^r u_3(s) q(s) ds
        \right\|_{\LebS{s(3,1)}} \\
    &\hspace{1cm} + \left\|
        - \frac{2n-2 + u_3 }{ r^2 } \nabla(  u_3(r) q(r) )ds
        \right\|_{\LebS{s(3,1)}} \\
    & =: A+B+C
\end{align*}

From the equation $u_{rr} + |u_r|^2 u = q_re$, we have
$|u_{rr}|\leq |q|^2 + |q|$ pointwise.
Therefore, for $A$,
\begin{align}
    A
    &\leq
    \left\|
        (|q|^2 + |q_r| ) \frac{1}{r^2} \int_0^r u_3(s)q(s)ds 
    \right\|_{\LebS{s(3,1)}} 
    \notag\\
    &=
    \| q \|_{\LebS{s(1,0)}} \left\| \frac{q}{r} \right\|_{\LebS{s(1,1)}}
    \left\|
        \frac{1}{r}\int_0^r u_3(s) q(s) ds
    \right\|_{\LebS{s(1,0)}}   +
    \| q_r \|_{\LebS{s(1,1)}} 
    \left\|
        \frac{1}{r^2}\int_0^r u_3(s) q(s) ds
    \right\|_{\LebS{s(2,0)}}
    \notag\\
    &=
    \| q \|_{\LebS{s(1,0)}} \|q/r\|_{\LebS{s(1,1)}}
    \left\|
        u_3(s) q(s) 
    \right\|_{\LebS{s(1,0)}}   +
    \| q_r \|_{\LebS{s(1,1)}} 
    \left\|
        \frac{1}{r}  u_3(s) q(s) ds
    \right\|_{\LebS{s(2,0)}}   
    \notag\\
    &\lesssim \|q\|_{\LebS{s(1,0)}}^2 \|q_r\|_{\LebS{s(1,1)}} + \| q_r \|_{s(1,1)} 
    \left\|
        \frac{1}{r} u_3(r) q(r)
    \right\|_{\LebS{s(2,0)}}   
    \notag \\
    &\lesssim \|q\|_{\LebS{s(1,0)}}^2 \|q_r\|_{\LebS{s(1,1)}} 
    + \| q_r \|_{s(1,1)} 
    \left\|
        \frac{u_3}{r} 
    \right\|_{\LebS{s(1,0)}}
    \left\|
         q(r)
    \right\|_{\LebS{s(1,0)}}
    \lesssim \| q_r \|_{\LebS{s(1,1)}}^3.
    \label{eqn:4:estA}
\end{align}

For $B$, we have,
\begin{align}
    B 
    &\lesssim \left\|
        \frac{2n-2 + u_3 }{ r }
        \right\|_{\LebS{s(1,0)}}
        \left\| \frac{1}{r^3} \int_0^r u_3(s) q(s) ds
        \right\|_{\LebS{s(2,1)}}
    \notag\\
    &\lesssim \left\|
        u_r
        \right\|_{\LebS{s(1,0)}}
        \left\| \frac{1}{r^2} u_3(r) q(r)
        \right\|_{\LebS{s(2,1)}}
    \notag\\
    &\lesssim \left\|
        u_r
        \right\|_{\LebS{s(1,0)}}
        \left\| \frac{u_3(r)}{r}  \right\|_{s(1,0)}
        \left\|
            \frac{q(r)}{r}
        \right\|_{\LebS{s(1,1)}}
    \notag\\
    &\lesssim \left\|
        u_r
        \right\|_{\LebS{s(1,0)}}^2
        \left\| q_r
        \right\|_{\LebS{s(1,1)}}
    \lesssim \| q_r \|_{\LebS{s(1,1)}}^3.
    \label{eqn:4:estB}
\end{align}

For $C$, we have,
\begin{align}
    C
    &\lesssim
    \left\|
         \frac{2n-2 + u_3 }{ r^2 } \nabla(  u_3(r) q(r) )ds
    \right\|_{\LebS{s(3,1)}}
    \notag\\
    & \lesssim
    \left\|
    \frac{2n-2 + u_3 }{ r }
    \right\|_{\LebS{s(1,0)}}
    \left(
    \left\|
        \frac{1}{r}\nabla(  u_3) q
    \right\|_{\LebS{s(2,1)}}  +
    \left\|
        \frac{1}{r}(  u_3) \nabla q
    \right\|_{\LebS{s(2,1)}}  
    \right) 
    \notag\\
    & \lesssim
    \left\|
        q
    \right\|_{\LebS{s(1,0)}}
    \left(
        \| \nabla u_3 \|_{\LebS{s(1,0)}} 
        \left\| \frac{q}{r} \right\|_{\LebS{s(1,1)}}
        + \left\| \frac{u_3}{r} \right\|_{\LebS{s(1,0)}}
        \|\nabla q \|_{\LebS{s(1,1)}}
    \right) 
    \notag\\
    &\lesssim \| q_r \|_{s(1,1)}^3.
    \label{eqn:4:estC}
\end{align}

The three estimates
\eqref{eqn:4:estA},
\eqref{eqn:4:estB} and
\eqref{eqn:4:estC}
together give the
estimate $\| \nabla N_1 \|_{L_x^{s(3,1)}} \lesssim \| q_r \|_{L_x^{s(1,1)}}^3$.

As for $N_2$, we have,
\begin{align}
    \| \nabla( \alpha q )\|_{\LebS{s(3,1)}}
    & \lesssim
    \| \alpha_r q \|_{\LebS{s(3,1)}}
    +
    \| \alpha q_r
    \|_{\LebS{s(3,1)}} 
    \notag\\
    &\lesssim
    \| \alpha_r\|_{\LebS{s(2,1)}} \| q \|_{\LebS{s(1,0)}}
    +
    \| \alpha \|_{\LebS{s(2,0)}} \| q_r \|_{\LebS{s(1,1)}}
    \notag
    \\
    &\lesssim \| \alpha_r \|_{\LebS{s(2,1)}} \| q_r \|_{\LebS{s(1,1)}}.
    \label{eqn:4:estN2.1}
\end{align}
Then, using the expression for $\alpha_r$ in \eqref{eqn:4:alphar} and
the fact that $u_3 \in L^\infty$,
\begin{align}
    \| \alpha_r \|_{s(2,1)} 
    &\lesssim \| q q_r \|_{\LebS{s(2,1)}}
    + \left\| \frac{ q^2 }{r} \right\|_{\LebS{s(2,1)}}
        +
    \left\|
         \frac{2n-2 + u_3 }{ r^2 } q \int_0^r u_3(s) q(s) ds
    \right\|_{s(2,1)}   
    \notag\\
    &\lesssim
    \| q \|_{\LebS{s(1,0)}}
    \left(
        \| q_r \|_{\LebS{s(1,1)}}
        + \left\| \frac{q}{r} \right\|_{\LebS{s(1,1)}}
    \right)+
         \left\| \frac{q}{r} \right\|_{\LebS{s(1,1)}}
    \left\|
         \frac{ 1 }{ r } \int_0^r u_3(s) q(s) ds
    \right\|_{\LebS{s(1,0)}}
    \notag\\
    &\lesssim
        \| q_r \|_{\LebS{s(1,1)}}^3
        +
    \left\| \frac{q}{r} \right\|_{\LebS{s(1,1)}}
    \left\|
        q
    \right\|_{\LebS{s(1,0)}}
    \lesssim
    \| q_r \|_{\LebS{s(1,1)}}^3.
    \label{eqn:4:estN2.2}
\end{align}
The estimates
    \eqref{eqn:4:estN2.1} and
    \eqref{eqn:4:estN2.2} give
 $\| \nabla N_2 \|_{L_x^{s(3,1)}} \lesssim \| q_r \|_{L_x^{s(1,1)}}^3$.
 and hence 
 \eqref{eqn:4:Nbounds}.
The estimate \eqref{eqn:4:Nbounds2} follows from an identical argument.
\end{proof}

Theorem \ref{thm:2} is thus established.

\section{The `real' heat flow case}
\label{sec:realheatflow}

In this section we will discuss what might be termed the `real' equivariant heat flow from $\Cn$ to $\CPn$.
In the case when $\alpha=1$ and $\beta=0$, that is, for the harmonic map heat flow, it is possible to make an ansatz which further reduces the problem.
In terms of the spherical coordinates,
\begin{align}
	u_t &= u_{rr} + u |u_r|^2 + \frac{2n-1}{r} u_r + \frac{ 2n-2+u_3 }{ r^2 } (e_3-u \langle u, e_3 \rangle	),	\label{eqn:6:pde} \\
	u(r,0) &= v(r) = (v_1(r),v_2(r),0) \in T_{e_3} \Sp^2 , \nonumber
\end{align}
this ansatz involves assuming that the initial data is valued in one great circle passing through the north pole;
	that is, the initial data is of of the form $c(r) e_3 + d(r) v_0$. (See Figure \ref{fig:greatcircle}.)
In this case for $t>0$ the solution will continue to be valued in the same great circle.
To see this, let $w = v \times e_3$
	and let $a(r,t) = u(r,t) \cdot w$.
By taking the inner product of equation \eqref{eqn:6:pde} with $a w$ we have
\begin{align*}
	a  a_t 	&= a a_{rr} - a^2 |u_r|^2 + \frac{2n-1}{r} a a_r + \frac{2n-2+u_3}{r^2} (- a^2 u_3 )	
	\leq a a_{rr}  + \frac{2n-1}{r} a  a_r +  \frac{2n-3}{r^2} a^2.
\end{align*}
By integrating this equation and using the Hardy inequality with best constant $4/d^2 = 4/(2n-2)^2$ we determine that,
\begin{align*}
	\frac{d}{dt} \frac{1}{2} \int_{\C^n} (a)^2 dx 	&\leq \sigma_{2n-1} \int_0^\infty a \der{}{r}( r^{2n-1} a_r ) dr + (2n-3) \left\| \frac{a}{r} \right\|^2_{L^2}  	\\
		& = - \sigma_{2n-1} \int_0^\infty (a_r)^2 r^{2n-1} dr + (2n-3) \frac{4}{ (2n-2)^2 } \| a_r \|^2_{L^2}	 \leq 0,
\end{align*}
and hence $a(r,t) = 0$ for all time.
The solution is therefore a linear combination of $v_0$ and $e_3$.

\newcommand{\thesphere}{}
\newcommand{\theradius}{}
\newcommand{\theangle}{}
\newcommand{\thephi}{}

\newcommand{\fmeridian}[2][]{

	\begin{scope}[rotate={90-atan(sin(\theangle)*tan(#2))}]
		
		\begin{scope}[yscale={sin(#2)*sqrt( (cos(#2)^2+sin(#2)^2 * sin(\theangle)^2)/(cos(#2)^2 + tan(\theangle)^2 ) )}]
			\draw[#1] (\theradius,0) arc(0:180:{\theradius}); 
		\end{scope}

	\end{scope}
}
\newcommand{\bmeridian}[2][]{

	\begin{scope}[rotate={90-atan(sin(\theangle)*tan(#2))}]
		
		\begin{scope}[yscale={sin(#2)*sqrt( (cos(#2)^2+sin(#2)^2 * sin(\theangle)^2)/(cos(#2)^2 + tan(\theangle)^2 ) )}]
			\draw[#1] (-\theradius,0) arc(180:360:{\theradius}); 
		\end{scope}

	\end{scope}
}

\newcommand{\fparallel}[2][]{	

	\begin{scope}[yscale=sin(\theangle)]


		\draw[#1,shift={(0,{\theradius*sin(#2)/tan(\theangle)})}] ({\theradius*cos(#2)*sin(acos(min(1,tan(\theangle)*tan(#2))))},{
				sin(\theangle)*cos(#2)*min(1,tan(\theangle)*tan(#2))*\theradius/sin(\theangle)})
			 arc({90-acos(min(tan(#2)*tan(\theangle),1))}:{-270+acos(min(tan(\theangle)*tan(#2),1))}:{cos(#2)*\theradius});

	\end{scope}
}

\newcommand{\bparallel}[2][]{	

	\begin{scope}[yscale=sin(\theangle)]

		\draw[#1] ({\theradius*cos(#2)*sin(acos(min(1,tan(\theangle)*tan(#2))))},{\theradius*cos(\theangle)*sin(#2)/sin(\theangle)+cos(#2)*min(1,tan(\theangle)*tan(#2))*\theradius})
			 arc({90-acos(min(tan(#2)*tan(\theangle),1))}:{90+acos(min(tan(\theangle)*tan(#2),1))}:{cos(#2)*\theradius});

	\end{scope}
}

\newenvironment{sphere}[2]{
	\renewcommand{\theradius}{#1}
	\renewcommand{\theangle}{#2}
	\renewcommand{\thesphere}[1]{\draw (0,0) circle(\theradius);}

	\coordinate (northpole) at (0,{cos(\theangle)*\theradius});
	\coordinate (southpole) at (0,{-cos(\theangle)*\theradius});
}{}

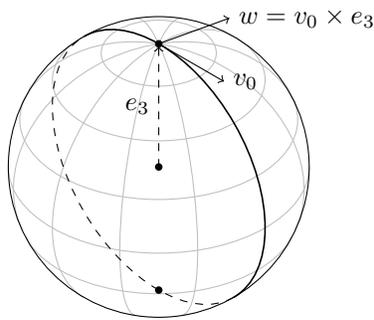
\begin{figure}
\centering
\begin{tikzpicture}

	\begin{sphere}{2}{35}

		\fmeridian[lightgray]{15}
		\fmeridian[lightgray]{-15}
		\fmeridian[lightgray]{45}
		\fmeridian[lightgray]{75}
		\fmeridian[lightgray]{-75}

		\fparallel[lightgray]{-22.5}
		\fparallel[lightgray]{0}
		\fparallel[lightgray]{22.5}
		\fparallel[lightgray]{45}
		\fparallel[lightgray]{67.5}

		\fmeridian[semithick]{-45}
		\bmeridian[dashed]{-45}

		\fill (0,0) circle(0.05);
		\fill (northpole) circle(0.05);
		\fill (southpole) circle(0.05);
	\draw[dashed,->] (0,0) -- ++(90:1.6) node[midway,anchor=east] {$e_3$};
	\draw[->] (northpole) -- ++(-30:1) node[anchor=west] {$v_0$};
	\draw[->] (northpole) -- ++(20:1) node[anchor=west] {$w=v_0\times e_3$};

		\thesphere{}
	\end{sphere}


\end{tikzpicture}

\caption[Visualization of the initial data for the harmonic map heat flow
    taking values on a fixed great circle.]
    {In the case of the harmonic map heat flow,
    if the initial data takes values in one great circle 
	(here the great circle spanned by $v_0$ and $e_3$),
        then the solution will be valued in the same great circle for future times.
    Both the harmonic maps and the self-similar solutions
    constructed in Section \ref{sec:selfsimilar} are of this type.}
	\label{fig:greatcircle}
\end{figure}

In terms of the stereographic representation of the problem,
\begin{align*}
	f_t 	&= f_{rr} - \frac{ 2 \of f_r^2 }{1+|f|^2}	+ \frac{2n-1}{r} f_r(r) + \frac{2n-1}{r^2} f(r) + \frac{ 2 |f|^2 f }{ 1+|f(r)|^2 },  	\\
	f(r,0) 	&= f_0(r),
\end{align*}
the ansatz is that the initial data is of the form $f(r,0) = b(r) e^{i\theta}$, for some real valued function $b(r)$ and a constant $\theta$.
The solution will then be
	of the form $f(r,t) = b(r,t) e^{i\theta}$, for the same constant $\theta$ and for some real valued function $b(r,t)$.
This motivates the terminology `real heat flow'.

It is not surprising that this problem is simpler to analyze, and in fact with this assumption we are able to say more about the dynamics of the problem.
On the other hand, this problem is still interesting because both the harmonic maps and the self-similar solutions constructed in Section \ref{sec:selfsimilar} fit into this context. 
In fact, the harmonic maps are given in the stereograpic coordinates by $f(r,t) = \alpha r = |\alpha| r e^{i\theta}$.
The initial data for a self-similar solution is just a point, so the initial data is valued in the great circle passing through that point and the north pole.

We will now describe how, based on the ansatz just described, a simpler PDE on the solution may be determined.
As the solution is valued on a great circle we can perform a change of variables,
$
	u(r,t) = \cos(g) e_3 + \sin(g) v_0,
$
for an unknown real-valued $g$.
Geometrically, $g$ is the spherical distance between $u(r,t)$ and $e_3$.
We calculate,
$
	u_r =  g_r ( -\sin(g) e_3 + \cos(g) v_0 ),
$
and,
\begin{align*}
	u_{rr} 	&= g_{rr} (-\sin(g) e_3 + \cos(g) v_0 ) + g_r^2 ( -\sin(g) e_3 - \cos(g) v_0 )	
        = g_{rr} (-\sin(g) e_3 + \cos(g) v_0 ) - u |u_r|^2.
\end{align*}
Substituting these into \eqref{eqn:6:pde} gives,
\begin{align*}
	g_t (-\sin(g) e_3 + \cos(g) v_0 )
 	&=\left( g_{rr} + \frac{2n-1}{r} g_r \right) 
		(-\sin(g) e_3 + \cos(g) v_0 )	\\
	&\hspace{1cm}
		+ \frac{ 2n-2+\cos(g) }{ r^2 } ( e_3 - \cos(g) ( \cos(g) e_3 + \sin(g) v_0 )).
\end{align*}
Taking the inner product of this equation with $-\sin(g)e_3 + \cos(g) v_0$ then yields the equation on $g$.
\begin{defn}
	The real heat flow problem is the Cauchy problem,
\begin{equation}
	\label{eqn:6:gpde}
	g_t = g_{rr} + \frac{2n-1}{r} g_r - \frac{1}{r^2} \left[ (2n-2) \sin(g) + \frac{1}{2} \sin(2g) \right],
\end{equation}
	subject to the initial condition $g(r,0)=g_0(r)$.
\end{defn}
For convenience we let $\eta(x) = (2n-2) \sin(x) + \sin(2x)/2$.

\begin{defn}
	The stationary real heat flow problem is the ODE,
	\begin{equation}
		\label{eqn:6:harmonic}
		0 =	\psi_\alpha''(r) + \frac{2n-1}{r} \psi_\alpha'(r) - \frac{1}{r^2} \eta(\psi_\alpha),
	\end{equation}
	subject the initial conditions $\psi_\alpha(0) = 0$ and $\psi_\alpha'(0) = \alpha >0	$.
\end{defn}
In the spherical coordinates the stationary solutions -- that is, the harmonic maps -- are given explicitly in \eqref{eqn:2:harmap}.
By transforming these solutions into the coordinates $g$, one finds that the unique solutions to the stationary real heat flow problem are,
\[
	\psi_\alpha(r) = 2 \arctan(\alpha r),
\]
which may be verified by substitution into \eqref{eqn:6:harmonic}.
In light of later results, what will be most  notable about the explicit solution is that it is independent of $n$.

\subsection{Uniqueness of solutions to the PDE problem in the $n\geq 3$ case}

PDEs of the type,
\begin{equation}
	\label{eqn:6:genpde}
	u_t = u_{rr} + \frac{d-1}{r} u_r - \frac{ \eta(u) }{ r^2 },
\end{equation}
with,
\begin{align*}
	\eta(0)	&=\eta(\pi)=\eta(2\pi),
	&\eta(x)	&>0\text{ for }x \in (0,\pi),
	&\eta(x)	&<0\text{ for }x \in (\pi,2\pi),
\end{align*}
arise naturally in the study of the equivariant harmonic map heat flow from $\R^d$ to spherically symmetric manifolds.
There is a general theorem classifying when there is uniqueness of solutions and when there is not uniqueness \cite{Germain2017}.
It states that if,
\begin{equation}
	\label{eqn:6:nonuniqcase}
	\eta'(\pi) <- \frac{ (d-2)^2 }{4},
\end{equation}
then there is non-uniqueness -- that is, two distinct solutions with the same initial data -- 
while if,
\begin{equation}
	\label{eqn:6:uniqcase}
	\eta'(x) \geq - \frac{ (d-2)^2 }{ 4 },
\end{equation}
for all $x$
then for every initial data there is at most one solution in $L^\infty_t L^\infty_x$.
We offer the following new proof of the latter case.

\begin{prop}
	Suppose that $\eta'(x) \geq -(d-2)^2/4$ for all $x$.
	There there is at most one solution to \eqref{eqn:6:genpde} in $L^\infty_tL^\infty_x$.
\end{prop}
\begin{proof}
First we observe that the condition \eqref{eqn:6:uniqcase} implies the one-sided Lipshitz inequality,
\[
	\frac{ \eta(u) - \eta(v) }{ u-v } \geq \min_{x \in [0,2\pi]} \eta'(x)  \geq - \frac{(d-2)^2}{4}.
\]

Now consider two solutions $u$ and $v$ of \eqref{eqn:6:genpde} with the same initial data $u_0$ and
	set $\phi=u-v$. 
We will assume that $u_0 \in L^2 \cap L^\infty$;
the argument to upgrade this to $L^2$ is standard \cite{Germain2017}.
Under this assumption we calculate,
\begin{align}
	\frac{1}{2} \frac{d}{dt} \| \phi \|_{L^2}^2
		&= \frac{1}{2} \frac{d}{dt} \sigma_{d-1} \int_0^\infty 
                    |\phi(r,t)|^2 r^{d-1} dr 	
		= \sigma_{d-1} \int_0^\infty \phi \left[ \phi_{rr} + \frac{2n-1}{r} \phi_r - \frac{ \eta(u)-\eta(v) }{r^2} \right] r^{2n-1}dr 	\nonumber	\\
		&= - \| \phi_r \|^2_{L^2} -  \sigma_{d-1} \int_0^\infty \frac{\phi^2}{r^2} \left[ \frac{ \eta(u)-\eta(v) }{u-v} \right] r^{2n-1}dr  
		\leq - \| \phi_r \|^2_{L^2} + \frac{ (d-2)^2 }{ 4 } 
                    \left\| \frac{ \phi }{ r } \right\|_{L^2}^2	\nonumber	\\
		&\leq - \| \phi_r \|^2_{L^2} + \frac{ (d-2)^2 }{ 4 } \frac{ 4 }{ (d-2)^2 } \left\|\phi_r  \right\|_{L^2}^2 \leq 0,
\end{align}
where in the last line we have used Hardy's inequality with the best constant $4/(d-2)^2$.
This implies that $\phi \equiv 0$, and hence that $u=v$.
\end{proof}

In this context of the real equivariant heat flow from $\Cn$ to $\CPn$, this implies the following result (given as Theorem 5 (i) in the introduction).
\begin{prop}
	Let $n\geq 3$. For a given initial data there is at most one solution to \eqref{eqn:6:gpde} in $L^\infty_t L^\infty_x$.
\end{prop}
\begin{proof}
	Here $d=2n$ and $\eta(x) = (2n-2) \sin(x) + \sin(2x)/2$.
	We calculate,
	\begin{align*}
		\eta'(x)	
				&= (2n-2) \cos(x) + 2\cos^2(x)-1	\\
				&= (2n-6) \cos(x) + 2(\cos(x)+1)^2 -3	
				\geq (2n-6)(-1) + 0 - 3 = -(2n-3),
	\end{align*}
	where the last inequality holds because $n\geq 3$ and so $(2n-6)\geq0$.
	Now using the inequality $-(2n-3) \geq -(n-1)^2$ (which is equivalent to $3 \geq - (n+1)^2$)  gives condition \eqref{eqn:6:uniqcase} and hence the result.
\end{proof}

\subsection{The $\CPn[2]$ case: breakdown of uniqueness}

\begin{figure}
	\centering
	\begin{tikzpicture}
	
		\draw[<->] (0,3.5) -- (0,-2.5);
		\draw[<->] (-0.5,0) -- (4.5,0);

		\draw[semithick,domain=0:4,samples=100] plot( {\x}, {2*cos(90*\x)+cos(2*90*\x)} );

		\draw[dashed] (0,-1) -- ++(4,0);
		\node[anchor=east] at (0,-1) {$-1$};

		\fill (2,-1) circle(0.05);
		\fill (2,0) circle(0.05) node[anchor=south] {$\pi$};
		\fill (0,0) circle(0.05) node[anchor=south west] {$0$};
		\fill (4,0) circle(0.05) node[anchor=south] {$2\pi$};

		\node[anchor=east] at (0,-2) {$-2$};
		\node[anchor=east] at (0,1) {$1$};
			\draw (0,1) -- ++(-0.1,0);
			\draw (0,2) -- ++(-0.1,0);
			\draw (0,3) -- ++(-0.1,0);
			\draw (0,-1) -- ++(-0.1,0);
			\draw (0,-2) -- ++(-0.1,0);
		\node[anchor=east] at (0,2) {$2$};
		\node[anchor=east] at (0,3) {$3$};

		\draw[<-] (3.75,-1.1) to[out=-60,in=60] ++(0,-0.9) node[anchor=east] {\footnotesize $-(n-1)^2$ threshold};

	\begin{scope}[shift={(6,0)}]
		\draw[<->] (0,3.5) -- (0,-2.5);
		\draw[<->] (-0.5,0) -- (4.5,0);

		\draw[semithick,domain=0:4,samples=100] plot( {\x}, {2*cos(90*\x)+(0.5)*cos(2*90*\x)} );

		\draw[dashed] (0,-2) -- ++(4,0);

		\fill (2,-1.5) circle(0.05);
		\fill (0,2.5) circle(0.05);
		\fill (4,2.5) circle(0.05);

		\foreach \n in {-4,-2,2,4,6} {
			\draw (0.07,\n/2) -- ++(-0.14,0) node[anchor=east] {$\n$};
		}

		\foreach \n in {2,4} {
			\draw (\n,0.07) -- ++(0,-0.14);
		}
		\node[anchor=south east] {$0$};
		\node[anchor=south] at (2,0.07) {$\pi$};
		\node[anchor=south] at (4,0.07) {$2\pi$};

		\node[anchor=north] at (2,-2) {\footnotesize $-(n-1)^2$ threshold};
	\end{scope}
	\end{tikzpicture}
\caption[Plots of the function $\eta'(x)$ in the case of the real equivariant
    heat flow from $\C^n$ to $\CPn$  in the cases $n=2$ and $n=3$.]
{Plots of the function $\eta'(x)$ in the case of the real equivariant heat flow from $\C^n$ to $\CPn$  in the cases $n=2$ (left) and $n=3$ (right).
	For the $n=3$ case, we easily see that $\eta$ satisfies the condition \eqref{eqn:6:uniqcase} with $d=2n$,
	and hence that uniqueness in $L^\infty_t L^\infty_x$ holds.
	For the $n=2$ case, we see that both \eqref{eqn:6:nonuniqcase}  and \eqref{eqn:6:uniqcase} 
		do not hold, so the case does not fit into the general classification theory.
}
	\label{fig:plots}

\end{figure}
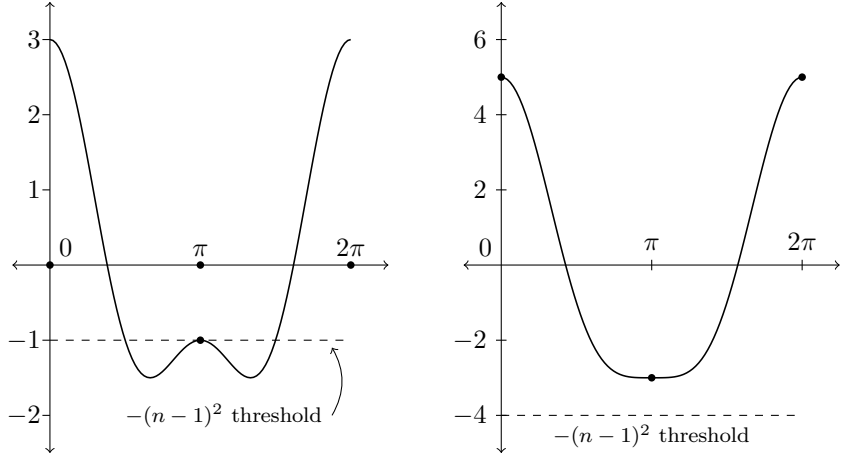

The $n=2$ case is the most interesting.
From the expression,
$
	\eta'(x) = 2\cos(x) + \cos(2x),
$
we see that $\eta'(\pi) = -1$, which is precisely the threshold $-(d-2)^2/4 = -1$ in the conditions \eqref{eqn:6:nonuniqcase} and \eqref{eqn:6:uniqcase}.
The condition that would imply non-uniqueness, \eqref{eqn:6:nonuniqcase}, does not hold.
However we find that,
\[
	\eta''(\pi) = - 2\sin(\pi) - 4 \sin(2\pi) = 0,
\]
and,
\[
	\eta'''(\pi) = -2\cos(\pi) - 8 \cos(2\pi) = -2(-1) -8(+1) = -6 <0,
\]
so in fact, by the second derivative test, $\pi$ is a local \emph{maximum} of $\eta'(x)$.
This means that the condition that would imply uniqueness, \eqref{eqn:6:uniqcase}, does hold either.
Hence the case of the real equivariant heat flow from $\C^2$ to $\CPn[2]$ is a borderline case 
	not covered by the classification theorem of \cite{Germain2011}.
(Plots of $\eta$ in the $n=2$ and $n=3$ cases are given in Figure \ref{fig:plots}, which make the difference clear.)

The question is then: does uniqueness hold or not?
First, we see that the proof of uniqueness presented in the last section clearly breaks down: because the derivative goes below the threshold value $-(d-2)^2/4$,
	a Lipshitz inequality of the form $\eta(u)-\eta(v)/(u-v) \geq -(d-2)^2/4$ cannot hold.

On the other hand, inspecting the proof in \cite{Germain2011} of non-uniqueness in the case \eqref{eqn:6:nonuniqcase} 
	we see that it relies critically on the following fact:
	if condition \eqref{eqn:6:nonuniqcase} holds, then the stationary solutions (that is, the harmonic maps) of the PDE problem oscillate around the fixed point $\pi$ as they converge to it.
In our case, the harmonic maps are given explicitly by $\psi_\alpha(r) = 2\arctan(\alpha r)$ and are clearly not oscillatory, and so that proof of non-uniqueness will not hold.
In fact, what is interesting is that the harmonic maps being monotonic is ordinarily a sign that there is uniqueness
	(if the uniqueness condition \eqref{eqn:6:uniqcase} holds, then the harmonic maps are necessarily monotonic.)
However, by using an alternative method in \cite{Germain2011} we are able to show that uniqueness for the problem from $\C^2$ to $\CPn[2]$ does not hold.
The original theorem requires some background to state, so we state a special version adapted to our setting.

\begin{thm*}[\cite{Germain2011}, Theorem 2.2]
	Suppose that the `equator map' $u(r,t) \equiv \pi$ (which is a time independent solution of the PDE)
	does not minimize the energy
	\[
		E(f) = \int_0^1 \left[ |f'|^2 + \frac{\gamma(f)}{r^2} \right] r^{d-1}dr
	\]
	where $\gamma'(x)=\eta(x)$.
	Then there exists a self-similar weak solution of the initial value problem \eqref{eqn:6:genpde} that is not constant in time and 
		that has the same initial data as the equator map, $u_0(r) \equiv \pi$.
\end{thm*}

Using this, we prove part (ii) of Theorem 5 in the introduction.

\begin{prop}
	For the case $n=2$ there is non-uniqueness of the problem \eqref{eqn:6:gpde}:
		there are two distinct solutions with initial data $u_0(r)\equiv\pi$.
\end{prop}

\begin{proof}
	The key aspect of the proof is capturing the fact that in the $n=2$ case, the condition $\eta'(x) \geq -(d-2)^2/4=-1$ in \eqref{eqn:6:uniqcase}
	is violated.
	If the non-uniqueness condition $\eta'(\pi)<-1$ in \eqref{eqn:6:nonuniqcase} held, this would be easy.
	However because $\eta'(\pi)=-1$, we need to do a higher order expansion of $\eta'(x)$ around $\pi$ to show this.
	Once we establish that condition $\eta'(x) \geq -1$ is violated, we follow \cite{Rupflin2010} and construct $h$
		based on a function which almost saturates that Hardy inequality.

	Let $u=\pi$ denote the equator map and $h$ be any function.
	We have
	\begin{equation}
		\label{eqn:6:energy}
		E(h) - E(u) =  \int_0^1 \left[ |h'|^2 + \frac{\gamma(h) - \gamma(\pi) }{r^2} \right] r^{d-1}dr,
	\end{equation}
	where $\gamma'(x) = \eta(x)$.
	One calculates,
	\begin{align*}
		\gamma'(\pi) 	&= \eta(\pi) = 0;	
		&\gamma''(\pi) 	&= \eta'(\pi) = -1;	
		&\gamma'''(\pi) 	&= \eta''(\pi) = 0; 
		&\gamma''''(\pi) 	&= \eta'''(\pi) = -6.
	\end{align*}
	Therefore by a Taylor expansion, if we choose $\delta$ small then there exists a constant $C>0$ such that	
	\begin{equation}
		\label{eqn:6:taylor}
		\gamma(x) - \gamma(\pi) \leq - (x-\pi)^2 - C(x-\pi)^4
	\end{equation}
	for all $x \in [\pi-\delta,\pi+\delta]$.
	The constant $C$ is positive because $\gamma^{(4)}(\pi)<0$.

	To use the inequality \eqref{eqn:6:taylor} in the energy expression \eqref{eqn:6:energy}, we need to choose $h$ valued in $[\pi-\delta,\pi+\delta]$.

	Following \cite{Rupflin2010}, we define, for any $\epsilon>0$, the function $f_\epsilon: [0,1] \rightarrow \R$ by
	\begin{equation}
		\label{eqn:6:fdefn}
		f_\epsilon(r) =
		\begin{cases}
			\epsilon^{-1}	&\text{for }0\leq r \leq \epsilon,	\\
			r^{-1}	&\text{for }\epsilon\leq r \leq 1/2,	\\
			4(1-r)	&\text{for }1/2\leq r \leq 1.
		\end{cases}.
	\end{equation}
	One verifies that $f(r)$ satisfies,
	\begin{equation}
		\label{eqn:6:hardyclose}
		\int_0^1 \left| \frac{f}{r} \right|^2 r^{3} dr \leq 
		\int_0^1 \left| {f'} \right|^2 r^{3} dr \leq 
		\left(1+ \frac{B}{|\log(\epsilon)|} \right) \int_0^1 \left| \frac{f}{r} \right|^2 r^{3} dr,
	\end{equation}
	for some $B>0$ independent of $\epsilon$. That is, $f$ is close to saturating the Hardy inequality, which in this case has best constant $4/(d-2)^2 = 1$.
	Then set	
	\[
		h(r) = \pi - \delta  \frac{ f_\epsilon(r) }{ \| f \|_{L^\infty} }
		 = \pi - \delta  \frac{ f_\epsilon(r) }{ 2 }.
	\]
	We observe that $h(r) \in [\pi-\delta,\pi+\delta]$ for all $r$.

	We then have
	\begin{align*}
		E(h) - E(u) 
			&=  \int_0^1 \left[ |h'|^2 + \frac{\gamma(h) - \gamma(\pi) }{r^2} \right] r^{d-1}dr	
		    \leq  \int_0^1 \left[ |h'|^2 + \frac{-(h-\pi)^2 - C (h-\pi)^4 }{r^2} \right] r^{d-1}dr	\\
			&=  \int_0^1 \left[  \frac{\delta^2}{4} |f_\epsilon'|^2 -  \frac{\delta^2}{4}  \left| \frac{f_\epsilon}{r} \right|^2 
					- C  \frac{\delta^4}{ 16 }   \left| \frac{f_\epsilon}{r} \right|^2 |f_\epsilon|^2 
 				\right] r^{d-1}dr	
	\end{align*}
	Now using bound \eqref{eqn:6:hardyclose} we determine that
	\begin{align*}
		E(h) - E(u) \leq
			&  \frac{\delta^2}{ 4 }  
			 \int_0^1 \left| \frac{f_\epsilon}{r} \right|^2 \left(    \frac{ B }{ |\log \epsilon| }    - C \frac{ \delta^2 }{ 4 } |f_\epsilon|^2 \right)  r^{3}dr,
	\end{align*}
	and by choosing $\epsilon$ sufficiently small we may make the right hand side negative.

	We thus determine that $E(h)<E(u)$, and hence there are two solutions.
\end{proof}

\subsection{The $n\geq 3$ case: precise  dynamics of the self similar solutions}

We finally present some results on the dynamics of the self-similar solutions in the real heat flow case when $n\geq 3$.
The methods of analysis here are not original, and our results are  based on analogous  results elsewhere.
Our motivation in presenting them here to show how in this special case, one can determine precise  dynamics of the self-similar solutions;
	it would be very satisfactory to extend these results to the general case of the \gll{} equation.

We first recall the self-similar problem.
\begin{defn}
	The self-similar real heat flow problem is the ODE
	\begin{align*}
		0 &=	\phi_\beta''(r) + \left( \frac{2n-1}{r} +\frac{r}{2} \right) \phi_\beta'(r) - \frac{1}{r^2} \eta(\phi_\beta)
	\end{align*}
	subject the initial conditions $\phi_\beta(0) = 0$ and $\phi_\beta'(0) = \beta >0	$.
\end{defn}

From section 2 we know that for every $\beta>0$ there is a unique global solution to this problem
and that there exists $\phi_\beta(\infty) \in \R$ such that $\lim_{r \rightarrow \infty} \phi_\beta(r) = \phi_\beta(\infty)$.

\begin{prop}
	\label{prop:precise}
	Let $\phi_\beta$ be the solution of the self-similar problem and $\psi_\beta$ the solution of the stationary problem.
	\begin{enumerate}
		\item[(i)] We have the bound $\phi_\beta(r) \leq \psi_\beta(r)$.

		\item[(ii)] The function $\phi_\beta$ is monotonically increasing and $\phi_{\beta}(r) < \pi$.

		\item[(iii)] For fixed $r>0$, the function $\beta \mapsto \phi_{\beta}(r)$ is strictly increasing, $\phi_0(r)=0$,
			and $\displaystyle \lim_{\beta \rightarrow \infty} \phi_\beta(r) = \pi$.

		\item[(iv)] The function $\beta \mapsto \phi_{\beta}(\infty)$ is strictly increasing, $\phi_0(\infty)=0$,
			and $\displaystyle \lim_{\beta \rightarrow \infty} \phi_\beta(\infty) = \pi$.

	\end{enumerate}
\end{prop}

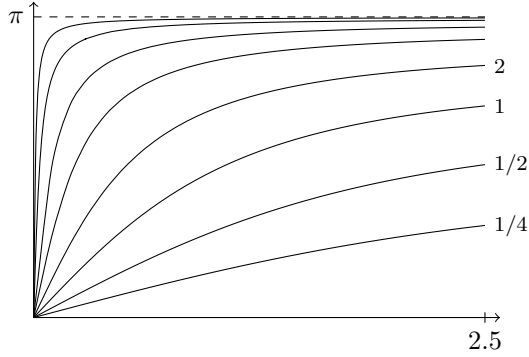
\begin{figure}
	\centering

	\begin{tikzpicture}


	\draw[ <-> ] (0,4.2) -- (0,0) -- (6.2,0);
	\draw[dashed] (0,4.0) -- ++(6,0);
	\node[anchor=east] at (0,4) {$\pi$};

	\draw (6,0.07) -- ++(0,-0.14) node[anchor=north] {$2.5$};

	\node[anchor=west] at (6,1.225) {\footnotesize ${1/4}$};
	\node[anchor=west] at (6,2.03) {\footnotesize ${1/2}$};
	\node[anchor=west] at (6,2.815) {\footnotesize ${1}$};
	\node[anchor=west] at (6,3.35) {\footnotesize ${2}$};

	\draw plot[smooth] coordinates{ (0,0)
(0.24, 0.0636288) (0.48, 0.127059) (0.72, 
  0.190096) (0.96, 0.252548) (1.2, 0.314234) (1.44, 
  0.374979) (1.68, 0.434624) (1.92, 0.49302) (2.16, 
  0.550035) (2.4, 0.605552) (2.64, 0.659469) (2.88, 
  0.711703) (3.12, 0.762184) (3.36, 0.81086) (3.6, 
  0.857695) (3.84, 0.902665) (4.08, 0.945763) (4.32, 
  0.986992) (4.56, 1.02637) (4.8, 1.06391) (5.04, 1.09966) (5.28, 
  1.13365) (5.52, 1.16594) (5.76, 1.19657) (6., 1.2256)
};

	\draw plot[smooth] coordinates{ (0,0)
(0.24, 0.127178) (0.48, 0.253488) (0.72, 
  0.378093) (0.96, 0.500209) (1.2, 0.619137) (1.44, 
  0.734269) (1.68, 0.845106) (1.92, 0.951259) (2.16, 
  1.05245) (2.4, 1.14849) (2.64, 1.23931) (2.88, 1.32488) (3.12, 
  1.40528) (3.36, 1.48062) (3.6, 1.55106) (3.84, 1.61679) (4.08, 
  1.67803) (4.32, 1.735) (4.56, 1.78795) (4.8, 1.83712) (5.04, 
  1.88274) (5.28, 1.92505) (5.52, 1.96429) (5.76, 2.00066) (6., 
  2.03437)
};
	\draw plot[smooth] coordinates{ (0,0)
(0.24, 0.253725) (0.48, 0.502047) (0.72, 
  0.740189) (0.96, 0.964452) (1.2, 1.17241) (1.44, 1.36288) (1.68,
   1.53565) (1.92, 1.69129) (2.16, 1.83082) (2.4, 1.95554) (2.64, 
  2.06683) (2.88, 2.16607) (3.12, 2.25458) (3.36, 2.33358) (3.6, 
  2.40417) (3.84, 2.46733) (4.08, 2.52394) (4.32, 2.57475) (4.56, 
  2.62046) (4.8, 2.66163) (5.04, 2.6988) (5.28, 2.73241) (5.52, 
  2.76287) (5.76, 2.79051) (6., 2.81565)
};
	\draw plot[smooth] coordinates{ (0,0)
(0.24, 0.502509) (0.48, 0.967821) (0.72, 1.37281) (0.96,
   1.71138) (1.2, 1.98862) (1.44, 2.21413) (1.68, 2.39787) (1.92, 
  2.54847) (2.16, 2.67294) (2.4, 2.77669) (2.64, 2.86394) (2.88, 
  2.93792) (3.12, 3.00113) (3.36, 3.05552) (3.6, 3.10263) (3.84, 
  3.14366) (4.08, 3.1796) (4.32, 3.21124) (4.56, 3.2392) (4.8, 
  3.26403) (5.04, 3.28616) (5.28, 3.30596) (5.52, 3.32372) (5.76, 
  3.33971) (6., 3.35414)
};
	\draw plot[smooth] coordinates{ (0,0)
(0.24, 1.07649) (0.48, 1.86432) (0.72, 2.37228) (0.96, 
  2.7015) (1.2, 2.92457) (1.44, 3.08292) (1.68, 3.19992) (1.92, 
  3.28924) (2.16, 3.35928) (2.4, 3.4154) (2.64, 3.46119) (2.88, 
  3.49911) (3.12, 3.53092) (3.36, 3.55789) (3.6, 3.58098) (3.84, 
  3.6009) (4.08, 3.61821) (4.32, 3.63334) (4.56, 3.64666) (4.8, 
  3.65843) (5.04, 3.66888) (5.28, 3.6782) (5.52, 3.68654) (5.76, 
  3.69403) (6., 3.70078)
};
	\draw plot[smooth] coordinates{ (0,0)
(0.24, 1.9996) (0.48, 2.81762) (0.72, 3.1774) (0.96, 
  3.37127) (1.2, 3.49081) (1.44, 3.57129) (1.68, 3.62889) (1.92, 
  3.67196) (2.16, 3.70527) (2.4, 3.73169) (2.64, 3.7531) (2.88, 
  3.77074) (3.12, 3.78548) (3.36, 3.79794) (3.6, 3.80857) (3.84, 
  3.81773) (4.08, 3.82568) (4.32, 3.83263) (4.56, 3.83873) (4.8, 
  3.84412) (5.04, 3.8489) (5.28, 3.85316) (5.52, 3.85697) (5.76, 
  3.86039) (6., 3.86348)
};
	\draw plot[smooth,tension=1] coordinates{ (0,0)
(0.024, 0.742467) (0.048, 1.3766) (0.072, 1.86656) (0.096, 
  2.23129) (0.12, 2.50302) (0.144, 2.70895) (0.168, 
  2.86848) (0.192, 2.99477) (0.216, 3.09677) (0.24, 
  3.18061) (0.264, 3.2506) (0.288, 3.30983) (0.312, 
  3.36054) (0.336, 3.40442) (0.36, 3.44274) (0.384, 
  3.47647) (0.408, 3.50639) (0.432, 3.53309) (0.456, 
  3.55706) (0.48, 3.5787) (0.504, 3.59833) (0.528, 
  3.61621) (0.552, 3.63256) (0.576, 3.64757) (0.6, 3.6614) (0.624,
   3.67418) (0.648, 3.68603) (0.672, 3.69703) (0.696, 
  3.70729) (0.72, 3.71686)
	(0.96, 
  3.78636) (1.2, 3.82801) (1.44, 3.85566) (1.68, 3.87529) (1.92, 
  3.88989) (2.16, 3.90114) (2.4, 3.91005) (2.64, 3.91726) (2.88, 
  3.92319) (3.12, 3.92814) (3.36, 3.93232) (3.6, 3.93589) (3.84, 
  3.93896) (4.08, 3.94163) (4.32, 3.94395) (4.56, 3.946) (4.8, 
  3.9478) (5.04, 3.94941) (5.28, 3.95083) (5.52, 3.95211) (5.76, 
  3.95326) (6., 3.95429)
};
	\draw plot[smooth,tension=0.8] coordinates{ (0,0)
	(0.024, 2.00568) (0.048, 2.82386) (0.072, 3.18404) (0.096, 
  3.37879) (0.12, 3.49945) (0.144, 3.58121) (0.168, 
  3.64015) (0.192, 3.68461) (0.216, 3.71933) (0.24, 
  3.74717) (0.264, 3.77) (0.288, 3.78904) (0.312, 3.80517) (0.336,
   3.81901) (0.36, 3.83101) (0.384, 3.84151) (0.408, 
  3.85078) (0.432, 3.85903) (0.456, 3.8664) (0.48, 
  3.87304) (0.504, 3.87904) (0.528, 3.8845) (0.552, 
  3.88949) (0.576, 3.89406) (0.6, 3.89826) (0.624, 
  3.90214) (0.648, 3.90573) (0.672, 3.90906) (0.696, 
  3.91216) (0.72, 3.91506)
  (0.84, 3.92704) (0.96,  3.93601) (1.08, 3.94297) (1.2, 3.94852)
	(1.44, 3.95681) (1.68, 3.96269) (1.92, 
  3.96706) (2.16, 3.97043) (2.4, 3.9731) (2.64, 3.97526) (2.88, 
  3.97703) (3.12, 3.97851) (3.36, 3.97976) (3.6, 3.98083) (3.84, 
  3.98175) (4.08, 3.98254) (4.32, 3.98324) (4.56, 3.98385) (4.8, 
  3.98439) (5.04, 3.98487) (5.28, 3.9853) (5.52, 3.98568) (5.76, 
  3.98602) (6., 3.98633)
};

	\end{tikzpicture}

\caption[Plots of the function $\phi_\beta(r)$ for $r \in [ 0, 2.5  ) $ and a variety of $\beta$ values.]
{Plots of $\phi_\beta(r)$ for $r \in [0,2.5]$ and $\beta=0.25$, $0.5$, $1$, $2$, $4.5$, $10$, $30$ and $100$.}
\label{fig:ss}
\end{figure}

The content of this Proposition may be seen at a glance in Figure \ref{fig:ss}.
Note that in light of the non-uniqueness theorem for $n=2$, we don't expect the same dynamics in the $n=2$ case:
	in fact we expect a self-similar profile whose asymptotic limit is $\pi$.

\newcommand{\pbe}{\psi_{\beta+\epsilon}}
\newcommand{\pb}{\phi_{\beta}}

\begin{lemma}
Suppose that for all $r\in [0,R]$, we have $\pb(r)<\pi$.
Then $\pb$ is increasing on $[0,R]$. 
\end{lemma}

\begin{proof}[Proof of lemma]
Becase $\beta>0$, the solution is initially increasing.
For a contradiction, let $r_0$ be the first critical  point in $[0,R]$.
Because $\pb$ is initially increasing, $r_0$ must be a local maximum.
However from the ODE we have
\[
	\pb''(r_0) = - \left( \frac{2n-1}{r} + \frac{r}{2} \right) \pb'(r_0) + \eta( \pb(r_0) ) = \eta(\pb(r_0)) >0,
\]
where $\eta( \pb(r_0) )>0$ because $\pb(r_0)\in(0,\pi)$.
The condition $\pb''(r_0)>0$ contradicts $r_0$ being a maximum.
 Hence $\pb$ is increasing on $[0,R]$.
\end{proof}

\begin{proof}[Proof of Proposition \ref{prop:precise}, (i)]
Let $\epsilon>0$ and consider the functions $\phi_\beta$ and $\psi_{\beta+\epsilon}(r)$.
Define
\[
	f(r) = r^2 ( \psi_{\beta+\epsilon}(r) - \phi_\beta(r) ).
\]
We will show that $f(r)\geq0$ for all $r$.
Letting $\epsilon \rightarrow 0$ will then give the result.

By continuity of derivatives given by the well-posedness theory, there is an initial interval $[0,\delta)$ on which $\psi_{\beta+\epsilon}(r) - \phi_\beta(r)$ is
	increasing, and hence, as $r^2$ is also increasing, the function $f$ is increasing on this interval.

Now suppose that $f$ has a critical point.
Let $r_0$ be the first critical point.
Because $f$ is initially increasing, this critical point must be a local maximum.
Because $f$ is increaing on $(0,r_0)$, we have $f(r_0)>0$.

We then calculate
\begin{align}
	f''(r) &= r^2 ( \psi_{\beta+\epsilon}''(r) - \phi_\beta''(r) )
			+ 4r ( \psi_{\beta+\epsilon}'(r) - \phi_\beta'(r) )
			+ 2 ( \psi_{\beta+\epsilon}(r) - \phi_\beta(r) )\nonumber	\\
		&= \frac{4-(2n-1)}{r} f' + \frac{ 2(2n-1) - 6 }{r^2 }f + \frac{r^3}{2} \phi_{\beta}' + \eta(\psi_{\beta+\epsilon}) - \eta( \phi_\beta). \label{eqn:6:f''}
\end{align}

Firstly, we have the Lipshitz bound 
\[
	\eta(\psi_{\beta+\epsilon}(r_0))-\eta(\phi_\beta(r_0)) \geq - (2n-3) (\pbe(r_0)-\pb(r_0)),
\]
	where we have used the fact that $f(r_0)=\pbe(r_0)-\pb(r_0)>0$ to multiply across by $\pbe(r_0) - \pb(r_0)$.

Secondly, because $f(r_0)>0$, $\pb(r_0) < \pbe(r_0) <\pi$, and hence by the Lemma $\pb$ is increasing on $[0,r_0]$.
Therefore $\pb'(r_0) \geq 0$.

Using both of these inequalities, and as well as  $f'(r_0)=0$, in \eqref{eqn:6:f''} yields
\begin{align*}
	f''(r_0)	&\geq + \frac{ 2(2n-1) - 6 }{r_0^2} f(r_0) + 0 - \frac{2n-3}{r_0^2} f(r_0)	\\
			&= \frac{2n-5}{r_0^2} f(r_0) > 0,
\end{align*}
which contradicts $r_0$ being a local maximum.
Hence $f$ has no critical points; it is increasing for all $r$.
In particular, it is always positive, so $\pb(r)<\pbe(r)$ for all $r$.
Taking the limit $\epsilon \rightarrow 0$ then gives  $\pb(r)\leq \psi_\beta(r)$.
\end{proof}

\begin{proof}[Proof of Proposition \ref{prop:precise}, (ii)]
The previous bound gives $\pb(r) \leq \psi_\beta(r) <\pi$ for all $r$.
Hence by the Lemma, $\pb(r)$ is always increasing.
\end{proof}

\begin{proof}[Proof of Proposition \ref{prop:precise}, (iii)]
Set $\alpha<\beta$.
We wish to show that $\phi_\alpha(r) < \phi_\beta(r)$,
	which follows from a maximum principle analysis of $g(r) = r^2(\phi_\beta(r) - \phi_\alpha(r)$.
The analysis is similar to the proof of item 2.
The function $g$ is is initially increasing.
If $r_0$ denotes the first critical point, which must be a maximum, one calculates
\begin{align*}
	g''(r_0) 	&=  \left[ \frac{ 4n-8 }{ r^2 } + 1 \right] g(r_0) + \frac{ \phi_\beta(r_0) - \phi_\alpha(r_0) }{ r^2 }\\
		 	&\geq \left[ \frac{ 4n-8 }{ r^2 } + 1 \right] g(r_0) -  \frac{ 2n-3  }{ r^2 }g(r_0) = \left[ \frac{2n-5}{r^2} +1 \right] g(r_0) \geq 0,
\end{align*}
a contradiction. Therefore $g$ is increasing for all $r$, and in particular is positive, and hence $\psi_\beta(r)>\psi_\alpha(r)$.
\end{proof}

\begin{proof}[Proof of Proposition \ref{prop:precise}, (iv)]
The proof follows from a similar maximum principle argument as in the previous proof to show that the function 
$
	h(r) = (r/(2+r))^2 (\psi_\beta(r) - \psi_\alpha(r))
$
is increasing.
One then has, for $r>1$,
\[
	\left( \frac{r}{2+r} \right)^2(\psi_\beta(r) - \psi_\alpha(r))
	\geq \frac{1}{9} (\psi_\beta(1) - \psi_\alpha(1)) >0,
\]
and hence on taking limits
$
	(\psi_\beta(\infty) - \psi_\alpha(\infty))
	\geq (1/9)(\psi_\beta(1) - \psi_\alpha(1)) >0,
$
which is what we wanted to prove.
\end{proof}

\appendix

\section{Some standard results}

\subsection{Hardy inequalities}

\begin{thm}[Generalized radial Hardy inequality]
    Suppose that $f: \R^{d} \rightarrow \R$
    is radial.
    Then for all $p \geq 1$ and $k \geq 0$ such that
    $p < d/(k+1)$ there holds,
    \begin{equation}
        \left\| \frac{f}{r^{k+1}} \right\|_{L^p}
        \leq \frac{ p }{ d- p(k+1) }
        \left\| \frac{f_r}{r^{k}} \right\|_{L^p}.
        \label{eqn:app:hardy1}
    \end{equation}
    \label{thm:app:hardy1}
\end{thm}

\begin{proof}
    We suppose that $f$ is smooth and compactly supported.
    The result for arbitrary $f$ then follows from a standard density argument.

    We have,
    \[
        \frac{d}{dr} \left( \frac{f}{r^{k}} \right)
        = - k \frac{f}{r^{k+1}} 
        + \frac{ f_r }{ r^k }.
    \]
    Multiplying this equation by $(f/r^{k+1})^{p-1} r^{d-1}$
    and integrating over $[0,\infty$ yields,
    \[
        \int_0^\infty 
        \frac{d}{dr} \left( \frac{f(r)}{r^{k}} \right)
        \left( \frac{f}{r^{k+1}} \right)^{p-1} r^{d-1}
        = - \frac{k}{ s(d) } \left\| \frac{f}{r^{k+1}} \right\|_{L^p}^p
        + \int_0^\infty \frac{ f_r }{ r^k } 
        \left( \frac{f}{r^{k+1}} \right)^{p-1} r^{d-1},
    \]
    where $s(d)$ is the measure of the unit sphere in $\R^d$.
    Now performing integration by parts on the term on the left we find,
    \begin{equation}
        \int_0^\infty 
        \frac{d}{dr} \left( \frac{f}{r^{k}} \right)
        \left( \frac{f}{r^{k+1}} \right)^{p-1} r^{d-1} 
        =
        -\int_0^\infty 
        \left( \frac{f}{r^{k}} \right) 
        \frac{d}{dr} \left[ 
        \left( \frac{f}{r^{k+1}} \right)^{p-1} r^{d-1} \right] dr
        + \left.\left[
            \left( \frac{f}{r^{k+1}} \right)^{p} r^d
        \right] \right|_{r=0}^{r=\infty}.
        \label{eqn:A:hardyproof}
    \end{equation}
    The boundary term corresponding to $r=\infty$ is 0 because
    $f$ is compactly supported.
    For the $r=0$ term we find,
    \[
        \lim_{r \rightarrow 0}
            \left( \frac{f}{r^{k+1}} \right)^{p} r^d
        =\lim_{r \rightarrow 0}
            f(r)^p r^{d - p(k+1)} = 0,
    \]
    if $d -p(k+1)>0$. We therefore have,
    \begin{align*}
        \int_0^\infty 
        \frac{d}{dr}& \left( \frac{f}{r^{k}} \right)
        \left( \frac{f}{r^{k+1}} \right)^{p-1} r^{d-1}  
        =
        -\int_0^\infty 
        \left( \frac{f}{r^{k}} \right) 
        \frac{d}{dr} \left[ 
        f(r)^{p-1} r^{d-1-(p-1)(k+1)}
        \right] dr
        \\
        &=
        -\int_0^\infty 
        \left( \frac{f}{r^{k}} \right) 
        \left[ 
        (d-1-(p-1)(k+1)) f(r)^{p-1} r^{d-2-(p-1)(k+1)}
        + \right.
        \\
        &\hspace{2cm} \left. +
        (p-1)f(r)^{p-2} f_r(r) r^{d-1-(p-1)(k+1)} 
        \right] dr
        \\
        &=
        -\frac{(d-p(k+1) + k )}{s(d)}
        \left\| \frac{f}{r^{k+1}} \right\|_{L^p}^p
        - (p-1) \int_0^\infty \frac{ f_r }{ r^k } 
        \left( \frac{f}{r^{k+1}} \right)^{p-1} r^{d-1}.
    \end{align*}
    Substituting this into \eqref{eqn:A:hardyproof} and combining terms we get,
    \[
        (d-p(k+1)) 
        \left\| \frac{f}{r^{k+1}} \right\|_{L^p}^p
        = -s(d) p \int_{\R^d} 
        \frac{ f_r }{ r^k } 
        \left( \frac{f}{r^{k+1}} \right)^{p-1} dx
        \leq p 
        \left\| \frac{f_r}{r^{k+1}} \right\|_{L^p}
        \left\| \frac{f}{r^{k+1}} \right\|_{L^p}^{p-1},
    \]
    which upon dividing through by the norm of $f/r^{k+1}$ gives
    the result.
\end{proof}

\begin{corr}
    Suppose that $f: \R^{d} \rightarrow X$
    is radial with $X=\C$ or $X=\R^m$.
    Then for all $p \geq 1$ and $k \geq 0$ such that
    $p < d/(k+1)$ there is a constant $C(d,p,X)$ such that,
    \begin{equation}
        \left\| \frac{f}{r^{k+1}} \right\|_{L^p}
        \leq C(d,p,X)
        \left\| \frac{f_r}{r^{k}} \right\|_{L^p}.
        \label{thm:app:hardy2}
    \end{equation}
\end{corr}

\begin{proof}
Take $X=\C$ and write $f$ as 
$f(r) = a(r) + i b(r)$ for real valued functions $a$ and $b$.
Using that
$\| u \|_{L^p} \sim \| \Re u \|_{L^p} + \| \Im u \|_{L^p}$,
we have,
\[
        \left\| \frac{f}{r^{k+1}} \right\|_{L^p}
        \lesssim
        \left\| \frac{a}{r^{k+1}} \right\|_{L^p} + 
        \left\| \frac{b}{r^{k+1}} \right\|_{L^p}  
        \lesssim  
        \left\| \frac{a_r}{r^{k}} \right\|_{L^p} +
        \left\| \frac{b_r}{r^{k}} \right\|_{L^p} 
        \lesssim
        \left\| \frac{f_r}{r^{k}} \right\|_{L^p}.
\]

A similar argument holds $X=\R^m$
writing $f$ in terms of its real-valued coordinate functions.
\end{proof}

\subsection{Local wellposedness for a class of singular ODE}

\begin{thm}
Consider the Cauchy problem,
\begin{align}
	f''(r) &= A\left(f'(r),f(r),r\right) - k
		\left( \frac{f'(r)}{r} - \frac{ f(r) }{ r^2 } \right) + \frac{1}{r^2} B(f(r)), \label{eqn:app:Cauchy}	\\
	f(0)	&=0,	\nonumber \\
	f'(0)	&=\alpha \in \C ,\nonumber
\end{align}
where
\begin{itemize}
	\item $k > 0$, 
	\item $A(z_1,z_2,r)$ is a smooth function with $A(\alpha,0,0)=0$,
	\item $B(z)$ is a smooth function such $|B(z)| \leq C |z|^3$ in a neighbourhood of $0$, and 
		$(\partial B/ \partial z)(0) = (\partial B/\partial \overline{z})(0)=0$.
\end{itemize}
There exists $r_0>0$ such that there is a unique solution among all functions $f:[0,r_0] \rightarrow \C$ satisfying,
\begin{equation}
	\left| f(r) \right|_{L^{\infty}{([0,r_0])}	}	+
	\left| \frac{f'(r) - f'(0)}{r} \right|_{L^{\infty}{([0,r_0]) } } < \infty.
	\label{eqn:app:condonf}
\end{equation}
The unique solution in this space is second differentiable at $r=0$ and satisfies $f''(0)=0$.
\label{thm:app:wellposedness}
\end{thm}

Let us make two remarks on the conditions in the theorem.
\begin{itemize}
	\item The condition \eqref{eqn:app:condonf} on $f$ is equivalent to both $f$ and $f'$ belonging to $L^\infty$ and $f'$ satisfying a  Liphitz condition at $r=0$.
	\item The assumptions on $B$ ensure that its behaviour as $r \rightarrow 0$ is non-singular; indeed, one readily verifies
		that, for smooth $f$, $B(f(r))/r^2 \rightarrow 0 $ as $r\rightarrow 0$.
		With this formulation of the Cauchy problem the singular behavior occurs only in the term $ \kappa (f'(r)/r - f(r)/r^2)$.
\end{itemize}

The proof the Theorem involves a standard, if delicate,
fixed point argument; details may be found in \cite{MyThesis}.

\subsection{An integration inequality}
\begin{prop}
	\label{prop:app:integration}
	Suppose that $A'(r) + c_1 r A(r) \leq c_2r^{-k}$ for $c_1>0$.
	Then for any $r_0>0$, $A(r) \leq C(c_1,r_0)(A(r_0)e^{-c_1r^2/4} + c_2 r^{-k+1})$.
\end{prop}
\begin{proof}
	We may write the equation as,	
	\[
		\frac{d}{dr} \left( e^{ c_1 r^2/2 } A(r) \right) \leq \frac{c_2}{r^k} e^{ c_1 r^2/2 },
	\]
	which on integration gives,
	\begin{align*}
		A(r) 	&\leq e^{c_1 \left( r_0^2-r^2 \right)/2} A(1) + c_2 e^{-c_1 r^2/2} \int_{r_0}^r \frac{1}{s^k} e^{c_1 s^2/2} ds	\\
		 	&= e^{c_1 \left(r_0^2-r^2\right)/2} A(1) + c_2 \frac{1}{r^{k-1}} \left( \frac{1}{ r^{-k+1} e^{c_1 r^2/2} } \int_{r_0}^r \frac{1}{s^k} e^{c_1 s^2/2} ds \right).
	\end{align*}
	To prove the result we show that the term in the brackets is bounded independently of $r$.
	This term is clearly a continuous function of $r$.
	Moreover, we have from the condition $c_1>0$,
	\[
		\lim_{r \rightarrow \infty} r^{-k+1} e^{ c_1 r^2/2 } = \infty	\;\;\text{ and }
		\lim_{r \rightarrow \infty}  \int_{r_0}^r \frac{1}{s^k} e^{c_1 s^2/2} ,
	\]
	which means, by L'Hopital's rule, that,
	\begin{align*}
		& \lim_{r \rightarrow \infty}  \left( \frac{1}{ r^{-k+1} e^{c_1 r^2/2} } \int_{r_0}^r \frac{1}{s^k} e^{c_1 s^2/2} ds \right)	\\
			&\hspace{0.5cm}= \lim_{r \rightarrow \infty} \left( \frac{1}{ (-k+1)r^{-k} e^{c_1 r^2/2} + c_1 r^{-k+2} e^{c_1 r^2/2} } \cdot\frac{1}{r^k} e^{c_1 r^2/2}  \right)=
			\lim_{r \rightarrow \infty} \frac{1}{ -k+1+c_1r^2} =0.
	\end{align*}
	We thus have for all $r\in[r_0,\infty)$,
	\[
		  \left( \frac{1}{ r^{-k+1} e^{c_1 r^2/2} } \int_{r_0}^r \frac{1}{s^k} e^{c_1 s^2/2} ds \right) \leq C(r_0,c_1),
	\]
	which completes the proof.
\end{proof}

\section*{Acknowledgments}

I wish to thank my doctoral advisor Pierre Germain
for suggesting the topic of the work,
providing critical help during the research process,
and for his unremitting patience as the work developed.
My thanks also to Chongchun Zeng for very helpful feedback,
    especially with regards to possible techniques for proving 
    wellposedness as in Section \ref{sec:wellposedness}.

\bibliographystyle{amsplain}
\bibliography{references}

\end{document}